\documentclass[twoside]{article}
\usepackage{graphicx,amssymb,mathrsfs,amsmath}
\renewcommand{\baselinestretch}{1.06}

\catcode`@=11
\long\def\@makefntext#1{\noindent #1}
\newskip\tabcentering \tabcentering=1000pt plus 1000pt minus 1000pt
\def\REF#1{\par\hangindent\parindent\indent\llap{#1\enspace}\ignorespaces} 
\def\MCH#1#2{\setbox0=\hbox{\raise#1\hbox{#2}}\smash{\box0}}

\def\@evenfoot{}\def\@oddfoot{}

\def\@evenhead{\hbox to\textwidth{\footnotesize\rm\thepage \hfill
{\it Tiexin Guo}}} 

\def\@oddhead{\hbox to \textwidth{\footnotesize{\it
A comprehensive connection between
the basic results and
 properties derived from two kinds of
  topologies } \hfill\thepage}}

\catcode`@=12

\def\bc{\begin{center}}
\def\ec{\end{center}}
\def\no{\noindent}
\def\hang{\hangindent\parindent}
\def\textindent#1{\indent\llap{\qquad #1\ \ \enspace}\ignorespaces}
\def\ref{\par\hang\textindent}

\begin{document}
\abovedisplayskip=6pt plus 1pt minus 1pt \belowdisplayskip=6pt
plus 1pt minus 1pt
\thispagestyle{empty} \vspace*{-1.0truecm} \noindent
\vskip 10mm

\bc{\Large\bf A comprehensive connection between
the basic results\\ and
 properties derived from two kinds of
  topologies\\ for a
random
locally convex module

\footnotetext{\footnotesize $^{\ast}$
E-mail address: txguo@buaa.edu.cn }
\footnotetext{\footnotesize Supported by NNSF No. 10871016}} \ec

\vskip 5mm
\bc{\bf Tiexin Guo$^{\ast}$}\\
{\small LMIB and School of Mathematics and Systems Science,
Beihang University, \\
Beijing 100191, P.R. China.}\ec

\begin{center}\bf Dedicated to Professor Berthold Schweizer on his 80th
birthday
\end{center}

\vskip 1 mm
\noindent\hrulefill\\
\noindent{\small {\small\bf Abstract}

The purpose of this paper is to make a comprehensive connection
between the basic results and properties derived from the two
kinds of topologies (namely the $(\varepsilon,\lambda)-$topology
introduced by the author and the stronger locally $L^{0}-$convex
topology recently introduced by Filipovi$\acute{c}$ et. al) for a
random locally convex module. First, we give an extremely simple
proof of the known Hahn-Banach extension theorem for
$L^{0}-$linear functions as well as its continuous variants. Then
we give the essential relations between the hyperplane separation
theorems in [Filipovi$\acute{c}$ et. al, J. Funct.
Anal.256(2009)3996--4029] and a basic strict separation theorem in
[Guo et. al, Nonlinear Anal. 71(2009)3794--3804]: in the process
obtain a useful and surprising fact that a random locally convex
module with the countable concatenation property must have the
same completeness under the two topologies! Based on the relation
between the two kinds of completeness, we go on to present the
central part of this paper: we prove that most of the previously
established deep results of random conjugate spaces of random
normed modules under the $(\varepsilon,\lambda)-$topology are still
valid under the locally $L^{0}-$convex topology, which
considerably enriches financial applications of random normed
modules.

\vspace{1mm}\baselineskip 12pt

\no{\small\bf Keywords} \ \ Random locally convex modules; The
$(\varepsilon,\lambda)-$topology; The locally $L^{0}-$convex topology;
Random conjugate spaces; Hahn-Banach extension theorems; Hyperplane
separation theorems; The countable concatenation property;
Completeness

\noindent\hrulefill\\

\noindent{\bf Contents} \\
\normalfont{
\begin{trivlist}
\item[1.] Introduction, preliminaries and an outline of main results.
 \begin{enumerate}
 \item[1.0.] Introduction
 \item[1.1.] The two kinds of topologies for a random locally convex
 module
 \item[1.2.] The random conjugate spaces of a random locally convex
 module under ${\cal T}_{\varepsilon,\lambda}$ and ${\cal T}_{c}$
 \item[1.3.] The countable concatenation closure of a set and
 hyperplane separation theorems
 \item[1.4.] Completeness
 \item[1.5.] The theory of random conjugate spaces of random normed
 modules under ${\cal T}_{c}$
 \end{enumerate}
\item[2.] Hahn-Banach extension theorems
 \begin{enumerate}
 \item[2.1.] Hahn-Banach extension theorems for random linear
 functionals
 \item[2.2.] Hahn-Banach extension theorems for $L^{0}-$linear
 functions
 \item[2.3.] Hahn-Banach extension theorems for continuous $L^{0}-$linear
 functions in random locally convex modules under the two kinds of
 topologies
 \end{enumerate}
\item[3.] Hyperplane separation theorems
 \begin{enumerate}
 \item[3.1.] The countable concatenation hull, the countable
 concatenation property and relations between hyperplane
 separation theorems currently available
 \item[3.2.] The hereditarily disjoint probability and more general
 forms of hyperplane separation theorems
 \item[3.3.] Closed $L^{0}-$convex subsets with the countable
 concatenation property
 \item[3.4.] The relation between the two kinds of completeness of a
 random locally convex module under the two kinds of topologies
 \end{enumerate}
\item[4.] The theory of random conjugate spaces of random normed
modules under locally $L^{0}-$convex topology
 \begin{enumerate}
 \item[4.1.] Riesz$^{\prime}$s representation theorems and the important
 connection between random conjugate spaces and classical conjugate
 spaces
 \item[4.2.] Random reflexivity and the James theorem
 \item[4.3.] Banach-Alaoglu theorem and
 Banach-Bourbaki-Kakutani-$\check{\textmd{S}}$mulian theorem
 \end{enumerate}
\item[5.] Some further remarks on the $(\varepsilon,\lambda)-$topology
and the locally $L^{0}-$convex topology
\item[Acknowledgments]
\item[References]
\end{trivlist}
} \baselineskip 15pt
\renewcommand{\baselinestretch}{1.02}
\parindent=16pt  \parskip=2mm
\rm\normalsize\rm
\newpage
\section*{1. Introduction, preliminaries and an outline of main results}
\subsection*{1.0. Introduction}
\quad\, Random metric theory originated from the theory of
probabilistic metric spaces, cf. [13,14,26]. The original
definition of a random normed space was presented in [26, p.240]
in which the random norm of a vector is a nonnegative random
variable. The new version of a random normed space was presented
in [11], in which the random norm of a vector is the equivalence
class of a nonnegative random variable. Based on the new version,
we further presented in [11] the elaborated definition of a random
normed module, which was originally introduced in [7]. Random
normed modules lead to the definitive notion of the random
conjugate space of a random normed space, cf. [11], and also make
the theory of random conjugate spaces of random normed modules
have obtained a systematic and deep development, cf.
[8,10,16,19,23]. Subsequently, random locally convex modules were
presented in [14] and deeply studied in [20,24](called random
seminormed modules in [14,20,24]). In particular, using the theory
of random conjugate spaces we recently established a basic strict
separation theorem in random locally convex modules in [22,
Theorem 3.1]. It should also be pointed out that random locally
convex modules, including their special case---random normed
modules, in all our previous papers, are endowed with a natural
topology, called the $(\varepsilon,\lambda)-$topology. The
terminology ``$(\varepsilon,\lambda)-$topology'' was first employed
by B.Schweizer and A.Sklar in 1961 in their work on topologizing a
probabilistic metric space, cf. [26]. The
$(\varepsilon,\lambda)-$topology is both useful and natural: for
example, the $(\varepsilon,\lambda)-$topology on the algebra of
equivalence classes of random variables on a probability space is
exactly the topology of convergence in probability, making the
algebra a topological algebra and a random locally convex module
endowed with its $(\varepsilon,\lambda)-$topology a topological
module over the above topological algebra.

Motivated by financial applications, Filipovi$\acute{c}$, Kupper
and Vogelpoth recently presented in [4] locally $L^{0}-$convex
modules, in particular the locally $L^{0}-$convex topology, establishing their hyperplane separation theorems, cf.
[4, Theorems 2.6 and 2.8], and some basic results on convex analysis
over the modules, cf. [4, Theorems 3.7 and 3.8]. Besides, they also rediscovered the notions of random locally convex modules and random normed modules [4, 25], in particular they
utilized their nice gauge function to show that the theory of
Hausdorff locally $L^{0}-$convex modules are equivalent to that of
random locally convex modules endowed with the locally
$L^{0}-$convex topology, cf. [4, Theorem 2.4]. The results in [4, 25]
enough exhibit the crucial importance of the locally
$L^{0}-$convex topology. When the topology is applied to a random
locally convex module we only need to notice that the algebra of
equivalence classes of random variables is only a topological ring
under its locally $L^{0}-$convex topology (since it is, in
general, so strong that it is no longer a linear topology, as
pointed out in [4, 25]) and similarly that a random locally convex
module endowed with the locally $L^{0}-$convex topology is only a
topological module over the topological ring. While the locally
$L^{0}-$convex topology is much stronger (or, finer) than the
$(\varepsilon,\lambda)-$topology, there are many attractive and
exciting relations between the basic results and important
properties derived from the two kinds of topologies for a random
locally convex module, for example, all the random locally convex
modules (in particular random normed modules) that have played
important roles in both financial applications and theoretic
studies have the same random conjugate spaces and completeness
under the two topologies.

The central purpose of this paper is to exhibit the essential
relations by making a comprehensive connection between our
previous work of a random locally convex module endowed with the
$(\varepsilon,\lambda)-$topology and Filipovic, Kupper and
Vogelpoth's basic work [4,25] of a locally $L^{0}-$convex module
(equivalently, a random locally convex module endowed with the
locally $L^{0}-$convex topology). First, we give an extremely
simple proof of the known Hahn-Banach extension theorem for
$L^{0}-$linear functions as well as its continuous variants. Then
we prove that our basic strict separation theorem implies
Filipovi$\acute{c}$, Kupper and Vogelpoth's hyperplane separation
theorem \uppercase\expandafter{\romannumeral 2} but is independent
of their hyperplane separation theorem
\uppercase\expandafter{\romannumeral 1}, also give a general
variant of either of their separation theorems, allowing a
separation with an arbitrary probability rather than the only
probability one, and in the process obtain a useful and surprising
fact that a random locally convex module with the countable
concatenation property must have the same completeness under the
two topologies! Based on the nice relation between the two kinds
of completeness, observing that a random normed module has the
same random conjugate space under the two kinds of topologies we
further present the central part of this paper: we prove that most
of the previously established deep results of random conjugate
spaces of random normed modules under the
$(\varepsilon,\lambda)-$topology are still valid under the locally
$L^{0}-$convex topology; at the same time, motivated by an important example
constructed by Filipovi$\acute{c}$, Kupper and Vogelpoth for
financial applications, in this part we may construct a
surprisingly wide class of random normed modules and give a
unified representation theorem of random conjugate spaces of them,
in particular we can also prove that the representation is
isometric, which generalizes and strengthens an important result
of Kupper and Vogelpoth's. In fact, seen from both our previous
work of theoretic researches and the Filipovic, Kupper and
Vogelpoth's recent work [4,25], the theory of random normed
modules together with their random conjugate spaces has been and
will be the most important part of the theory of random locally
convex modules together with their applications to conditional
risk measures, thus Section 4 of this paper is, without doubt, the
most important part of this paper, since this section has given
most of the important and deep results of the theory of random
conjugate spaces of random normed modules under the locally
$L^{0}-$convex topology.  Finally, the principal results of this
paper enough convince people that the two kinds of topologies
should be simultaneously considered in the future development of
random locally convex modules together with their financial
applications.

The remainder of this paper is organized as follows: In Subsections
1.1 and 1.2 of Section 1, as preliminaries we recapitulate some
basic facts on the two kinds of topologies and random conjugate
spaces of random locally convex modules, respectively; in the other
subsections of Section 1 we give an outline of the main results of
this paper; following Section 1 we state and prove our main results
in Sections 2 to 4 of this paper according to the order of contents.
Finally in Section 5 we conclude this paper with some further
remarks explaining the reason for which the two kinds of topologies should be
simultaneously considered in the future study of random locally
convex modules.

\subsection*{1.1. The two kinds of topologies for a random locally convex
 module}
\quad\, To introduce the two kinds of topologies, let us
 recapitulate the related terminology and
notation.

Throughout this paper, $(\Omega,{\cal F},P)$ denotes a probability
space, $K$ the scalar field $R$ of real numbers or $C$ of complex
numbers, $\bar{R}=[-\infty,+\infty]$, $L^{0}({\cal F},\bar{R})$
the set of equivalence classes of extended real-valued ${\cal
F}-$measurable random variables on $\Omega$, $L^{0}({\cal F},K)$
the algebra of equivalence classes of $K-$valued ${\cal
F}-$measurable random variables on $\Omega$ under the ordinary
scalar multiplication, addition and multiplication operations on
equivalence classes.

It is well known from [3] that $L^{0}({\cal F},\bar{R})$ is a
complete lattice under the ordering $\leqslant$: $\xi\leqslant \eta$
iff $\xi^{0}(\omega)\leqslant\eta^{0}(\omega)$, for $P-$almost all
$\omega$ in $\Omega$ (briefly,a.s.), where $\xi^{0}$ and $\eta^{0}$
are arbitrarily chosen representatives of $\xi$ and $\eta$,
respectively. Furthermore, every subset A of $L^{0}({\cal
F},\bar{R})$ has a supremum, denoted by $\vee A$, and an
infimum, denoted by $\wedge A$, and there exist a sequence
$\{a_{n}\,|\,n\in N\}$ and a sequence ${\{b_{n}\,|\,n\in N}\}$ in $A$ such that
$\vee_{n\geqslant1}$ $a_{n}=\vee A$ and $\wedge_{n\geqslant1}$
$b_{n}=\wedge A$, where $N$ denotes the set of positive integers.
If, in addition, $A$ is directed upwards(downwards), then the above
$\{a_{n}\}$ ($\{b_{n}\}$) can be chosen as nondecreasing
(nonincreasing). Finally $L^{0}({\cal F},R)$, as a sublattice of
$L^{0}({\cal F},\bar{R})$, is complete in the sense that every
subset with an upper bound has a supremum.

Besides, throughout this paper we distinguish random variables
from their equivalence classes by means of symbols: for example,
$I_{A}$ denotes the characteristic function of the ${\cal F}-$
measurable set A, then we use $\tilde{I}_{A}$ for its equivalence
class. It is necessary when we apply the theory of lifting
property to random normed modules as in [13, 14] and apply the
theory of random normed modules to the theory of random operators,
cf.[14,18]. Therefore, for a set A of random variables, esssup(A)
and essinf(A) denote its respective essential supremum and
infimum, they are still random variables, we also reserve
esssup(${\cal E}$) and essinf(${\cal E}$) for the essential
supremum and infimum of a subfamily ${\cal E}$ of ${\cal F}$,
respectively, as in [4].

Specially, $L^{0}_{+}=\{\xi\in L^{0}({\cal
F},R)\,|\,\xi\geqslant0\}$, $L^{0}_{++}=\{\xi\in L^{0}({\cal
F},R)\,|\,\xi>0 ~on~ \Omega\}$, where for A $\in {\cal F}$,
``$\xi>\eta ~on~ A$" means $\xi^{0}(\omega)>\eta^{0}(\omega)$ a.s.
on A for any chosen representatives $\xi^{0}$ and $\eta^{0}$ of
$\xi$ and $\eta$, respectively. As usual, $\xi>\eta$ means
$\xi\geqslant\eta$ and $\xi\neq\eta$.

Notice that a $K-$valued $P-$measurable function is exactly an
$\hat{{\cal F}}-$measurable one, where $\hat{{\cal F}}$ denotes the
completion of ${\cal F}$ with respect to $P$, then one can easily
see that the symbol $L(P,K)$ in [11, 14, 19, 20, 22, 23, 24] amounts
to $L^{0}(\hat{{\cal F}},K)$, which is essentially identified with
$L^{0}({\cal F},K)$ as a set of equivalence classes. Besides, an
${\cal F}-$measurable function must be an $\hat{{\cal
F}}-$measurable one, and thus the following Definition 1.1 was
employed in [11, 14, 19, 20, 22, 23, 24] in a slightly general way.

In the following Definition 1.1, we adopt the terminologies
``$L^{0}-$seminorms and $L^{0}-$norms '' (they were defined as
``Module-absolutely homogeneous random seminorms and random
norms'', respectively, in our papers [11, 13, 14, 19, 20, 22]) and
notation ``$\|\cdot\|$'' from [4,25] for simplicity, but the
essence of Definition 1.1 is the same as the original one used in
[11,14,20,22,24].

{\noindent\bf Definition 1.1.} An ordered pair $(E,{\cal P})$ is
called {\it{a random locally convex module over $K$ with base
$(\Omega,{\cal F},P)$}} if $E$ is a left module over the algebra
$L^{0}({\cal F},K)$ and ${\cal P}$ is a family of mappings from $E$
to $L^{0}_{+}$ such that the following three axioms are satisfied:

 \hspace{2mm}(\lowercase \expandafter
{\romannumeral 1}) $\vee\{\|x\|~|~\|\cdot\|\in{\cal P}\}=0$ iff
$x=\theta$(the null element of $E$);

\hspace{1mm}(\lowercase \expandafter {\romannumeral 2}) $\|\xi\cdot
x\|=|\xi|\cdot\|x\|,\forall\xi\in L^{0}({\cal F},K), ~x\in E$ and
$\|\cdot\|\in\cal P$;

(\lowercase \expandafter {\romannumeral 3})
$\|x+y\|\leqslant\|x\|+\|y\|,\forall x, y\in E$ and
$\|\cdot\|\in\cal P$.

 Furthermore, a mapping $\|\cdot\|: E\rightarrow L^{0}_+$
 satisfying (\lowercase \expandafter {\romannumeral 2}) and (\lowercase \expandafter {\romannumeral
 3}) is called {\it an $L^{0}-seminorm$}; in addition, if $\|x\|=0$ also implies $x=\theta$,
 then it is called {\it{an $L^{0}-$norm}}, in which case ($E,\|\cdot\|$) is
 called {\it{a random normed module (briefly, an $RN$ module) over $K$ with
 base $(\Omega,{\cal F},P)$}}, which was first introduced in [11],
 and is a special case of a random locally convex module when $\cal
 P$ consists of the only one $L^{0}-$norm $\|\cdot\|$.

Before giving the following Remark 1.2, let us first mention the
notion of an $L^{0}-$normed module, which was introduced in [25,
Definition 2.1] and employed in [4]. An ordered pair $(E,\|\cdot\|)$
is called {\it an $L^{0}-$normed module} if $E$ is a left module
over the ring $L^{0}({\cal F},R)$ and $\|\cdot\|$ is a mapping from
$E$ to $L^{0}_{+}$ such that the following three axioms are
satisfied:

(1)\, $\|x\|=0$ iff $x=\theta$(the null element of $E$),

(2)\, $\|\xi x\|=|\xi|\|x\|,\forall\xi\in L^{0}({\cal F},R)$ and
$x\in E$,

(3)\, $\|x+y\|\leq \|x\|+\|y\|,\forall x,y\in E$,\\
where $\|x\|$ is called the $L^{0}-$ norm of $x\in E$.

Since $L^{0}({\cal F},R)$ is a linear space with the linear space
structure as usual and is a ring with the unit element, as the
following Remark 1.2 shows, a left module over the ring $L^{0}({\cal
F},R)$ is, naturally, a left module over the algebra $L^{0}({\cal
F},R)$. Consequently the notion of an $L^{0}-$normed module is
essentially equivalent to that of an $RN$ module over $R$ with base
$(\Omega,{\cal F}, P)$.

It should also be pointed out that as Filipovi\'{c}, Kupper and
Vogelpoth emphasized in [4,25], the algebra $L^{0}({\cal F},R)$ is
only a topological ring under the locally $L^{0}-$convex topology
introduced by them in [4](see also Proposition 1.5 of this paper),
since the locally $L^{0}-$convex topology is too strong to be
necessarily a linear topology on $L^{0}({\cal F},R)$. The fact only
means that the locally $L^{0}-$convex topology on the algebra
$L^{0}({\cal F},R)$ is compatible with the ring structure but not
compatible with the scalar multiplication operation of $L^{0}({\cal
F},R)$, but this does not at all hinder us from endowing random
locally convex modules and $RN$ modules with the locally
$L^{0}-$convex topology; on the contrary, the results obtained in
this paper for random locally convex modules (in particular for $RN$
modules) not only are still suitable for locally $L^{0}-$convex
modules (correspondingly, for $L^{0}-$normed modules) but also make
our previous work and the work in [4,25] naturally connected
together.

 {\noindent\bf Remark 1.2.}  In Definition 1.1, $E$ is, of course, a linear space
 over $K$, and the module multiplication is a natural extension of
 the scalar multiplication: $\alpha\cdot x=(\alpha\cdot 1)\cdot
 x,~\forall\alpha\in K$ and $x\in E$, where 1 is the unit element of
 $L^{0}({\cal F},K)$. Conversely, if $E$ is only a left module
 over the ring $L^{0}({\cal F},K)$, then $E$ is again a left
 module over the algebra $L^{0}({\cal F},K)$ if the scalar
 multiplication is defined by $\alpha\cdot x=(\alpha\cdot 1)\cdot x,\forall\alpha\in K$
and $x\in E$. Thus the notion of an $L^{0}-$normed module in [4, 25]
is equivalent to that of an $RN$ module. One of advantages of our
formulation of an $RN$ module will be reflected in Proposition 1.4
below when an $RN$ module is endowed with the
$(\varepsilon,\lambda)-$topology, and this formulation means that it is
an $RN$ space, and thus has more advantages: for example, we can
often convert a problem of an $RN$ space to one of an $RN$ module,
cf. [9,16], in particular, cf. [16, Lemma 3.2], and the theory of an
$RN$ module can be applied to the theory of random linear operators
and functional analysis in a direct and convenient fashion, cf.
[7,18,12,15].

{\noindent\bf Example 1.3.}  $L^{0}({\cal F},K)$ is an $RN$ module
over $K$ with base $(\Omega,{\cal F},P)$ if it is endowed with the
$L^{0}-$norm $\|x\|=|x|,\forall x\in L^{0}({\cal F},K)$. It is well
known that $L^{0}({\cal F},K)$ is a topological algebra over $K$
endowed with the topology of convergence in probability $P$, a local
base at $\theta$ of which is
$\{N_{\theta}(\varepsilon,\lambda)|~\varepsilon>0,0<\lambda<1\}$,
where $N_{\theta}(\varepsilon,\lambda)=\{x\in L^{0}({\cal
F},K)|~P\{\omega\in\Omega| |x|(\omega)<\varepsilon\}>1-\lambda\}$.

Since B. Schweizer and A. Sklar introduced the
$(\varepsilon,\lambda)-$topology into more abstract
spaces---probabilistic metric spaces, namely they introduced ``the
topology of convergence in probability on probabilistic metric
spaces'', in 1961, their idea is also suitable on many other
occasions, cf.[26]. The following Proposition 1.4 is exactly a
copy of their idea in the case of a random locally convex module.

From now on, for a random locally convex module $(E,{\cal P})$ and
 for any finite subfamily ${\cal Q}\subset {\cal P}$,
 $\|\cdot\|_{{\cal Q}}$ always denotes the $L^{0}-$seminorm defined by $\|x\|_{{\cal Q}}=\vee\{\|x\| \,|\, \|\cdot\|\in {\cal Q}\},\forall x\in
 E$, unless otherwise stated.

{\noindent\bf Proposition 1.4 ([14, 20, 24, 22]).} Let $(E,\cal
P)$ be a random locally convex module over $K$ with base
$(\Omega,{\cal F},P)$. For any $\varepsilon>0$, $0<\lambda<1$ and any finite subfamily $Q$ of ${\cal P}$, let
$N_{\theta}(Q,\varepsilon,\lambda)=\{x\in
E\,|\,P\{\omega\in\Omega~|~\|x\|_{Q}(\omega)<\varepsilon\}>1-\lambda\}$
and ${\cal
U}_{\theta}=\{N_{\theta}(Q,\varepsilon,\lambda)~|~Q\subset {\cal P}$
finite, $\varepsilon>0,0<\lambda<1$\}, then $\mathcal
{U}_{\theta}$ is a local base at $\theta$ of some Hausdorff linear
topology, called the {\it {$(\varepsilon,\lambda)-$topology}}
induced by $\cal P$. Further, we have the following statements:

(1) $L^{0}({\cal F},K)$ is a topological algebra over $K$ endowed
with its $(\varepsilon,\lambda)-$topology, which is exactly the
topology of convergence in probability $P$;

(2) $E$ is a topological module over the topological algebra
$L^{0}({\cal F},K)$ when $E$ and $L^{0}({\cal F},K)$ are endowed
with their respective $(\varepsilon,\lambda)-$topologies;

(3) A net $\{x_{\alpha},~ \alpha\in\wedge\}$ in $E$ converges in the
$(\varepsilon,\lambda)-$topology to $x\in E$ iff
$\{\|x_{\alpha}-x\|,~\alpha\in\wedge\}$ converges in probability $P$
to 0 for each $\|\cdot\|\in\cal P$.

From now on, for all random locally convex modules, their
$(\varepsilon,\lambda)-$topologies are denoted by ${\cal
T}_{\varepsilon,\lambda}$ unless there is a danger of confusion.

{\noindent\bf Proposition 1.5 ([4])}. Let $(E,\cal P)$ be a random
locally convex module over $K$ with base $(\Omega,{\cal F},P)$. For
any $\varepsilon\in L^{0}_{++}$ and $Q\subset \mathcal {P}$ finite,
let
$B_{Q}(\varepsilon)=\{x\in E~|~\|x\|_{Q}\leqslant\varepsilon\}$ and
${\cal U}_{\theta}=\{B_{Q}(\varepsilon)\,|\,Q\subset\cal P$ finite,
$\varepsilon\in L^{0}_{++}\}$. A set $G\subset E$ is called ${\cal
T}_c-$ open if for every $x\in G$ there exists some
$B_{Q}(\varepsilon)\in {\cal U}_{\theta}$ such that
$x+B_{Q}(\varepsilon)\subset G$. Let ${\cal T}_c$ be the family of
${\cal T}_c-$open subsets, then ${\cal T}_c$ is a Hausdorff topology
on $E$, called the locally $L^{0}-$convex topology induced by $\cal
P$. Further, the following statements are true:

(1) $L^{0}({\cal F},K)$ is a topological ring endowed with its
locally $L^{0}-$convex topology;

(2) $E$ is a topological module over the topological ring
$L^{0}({\cal F},K)$ when $E$ and $L^{0}({\cal F},K)$ are endowed
with their respective locally $L^{0}-$convex topologies;

(3) A net $\{x_{\alpha},~ \alpha\in\wedge\}$ in $E$ converges in the
locally $L^{0}-$convex topology to $x\in E$ iff
$\{\|x_{\alpha}-x\|,~ \alpha\in\wedge\}$ converges in the locally
$L^{0}-$convex topology of $L^{0}({\cal F},K)$ to $\theta$ for each
$\|\cdot\|\in\cal P$.

From now on, for all random locally convex modules, their locally
$L^{0}-$convex topologies are denoted by ${\cal T}_{c}$ unless there
is a possible confusion.

 ${\cal T}_{c}$ is called locally
$L^{0}-$convex because it has a striking local base ${\cal
U}_{\theta}=\{B_{Q}(\varepsilon)~|~Q\subset {\cal P}$ finite and
$\varepsilon\in L^{0}_{++}\}$, each member $U$ of which satisfies
the following:

\hspace{2mm}(\lowercase \expandafter {\romannumeral 1})
$L^{0}-$convex: $\xi\cdot x+(1-\xi)\cdot y\in U$ for any $x,y\in
U$ and $\xi\in L^{0}_+$ such that $0\leqslant\xi\leqslant1$;

\hspace{1mm}(\lowercase \expandafter {\romannumeral 2})
$L^{0}-$absorbent: there is $\xi\in L^{0}_{++}$ for each $x\in E$
such that $x\in\xi\cdot U$;

(\lowercase \expandafter {\romannumeral 3}) $L^{0}-$balanced:
$\xi\cdot x\in U$ for any $x\in U$ and any $\xi\in L^{0}({\cal
F},K)$ such that $|\xi|\leqslant1$.

Such an ordered pair $(E,{\cal T}_c)$ such that ${\cal T}_c$
possesses the above properties (\lowercase \expandafter
{\romannumeral 1}), (\lowercase \expandafter {\romannumeral 2}) and
(\lowercase \expandafter {\romannumeral 3}) is called a locally
$L^{0}-$convex module (see [4, Definition 2.2]); conversely, for
every locally $L^{0}-$convex module $(E,{\cal T})$, ${\cal T}$ can
also be induced by a family of $L^{0}-$seminorms on $E$ as above,
see [4, Theorem 2.4] for its proof. Thus the theory of Hausdorff
locally $L^{0}-$convex modules amounts to the theory of random
locally convex modules endowed with the locally $L^{0}-$convex
topology.

\subsection*{1.2. The random conjugate spaces of a random locally convex module under
${\cal T}_{\varepsilon,\lambda}$ and ${\cal T}_{c}$.}
\quad\, Given a random locally convex module ($E,\cal P$) over $K$ with
base ($\Omega,{\cal F},P$), generally ${\cal T}_{c}$ is much
stronger than ${\cal T}_{\varepsilon,\lambda}$. The random
conjugate spaces of ($E,\cal P$) under ${\cal
T}_{\varepsilon,\lambda}$ and ${\cal T}_{c}$, however, coincide if
$\cal P$ has the countable concatenation property, in particular
an $RN$ module has the same random conjugate space under ${\cal
T}_{\varepsilon,\lambda}$ and ${\cal T}_{c}$. To see this, let
$E^{\ast}_{\varepsilon,\lambda}=\{f:E\rightarrow L^{0}({\cal
F},K)\,| f$ is a continuous module homomorphism from $(E,{\cal
T}_{\varepsilon,\lambda})$ to $(L^{0}({\cal F},K),{\cal
T}_{\varepsilon,\lambda})\}$ and $E^{\ast}_{c}=\{f:E\rightarrow
L^{0}({\cal F},K)\,| f$ is a continuous module homomorphism from
$(E,{\cal T}_{c})$ to $(L^{0}({\cal F},K),{\cal T}_{c})\}$, which
are called the random conjugate spaces of $(E,\cal P)$ under
${\cal T}_{\varepsilon,\lambda}$ and ${\cal T}_{c}$, respectively.
It should be noticed that $E^{\ast}_c$ was already used in [4] in
anonymous way.

Let us recall from [4] that $\cal P$ is called having the countable
concatenation property if $\sum_{n\geq 1}\tilde
{I}_{A_{n}}\cdot\|\cdot\|_{Q_{n}}$ still belongs to $\cal P$ for any
countable partition $\{A_{n}\,|\,n\in N\}$ of $\Omega$ to ${\cal F}$
and any sequence $\{Q_{n}~|~ n\in N\}$ of finite subfamilies of
$\cal P$.

It is easy to see that a random linear functional $f:E\rightarrow L^{0}({\cal F},K)\in E^{\ast}_{c}$ iff there are $\xi\in
L^{0}_+$ and $Q\subset \mathcal {P}$ finite such that
$|f(x)|\leqslant\xi\cdot\|x\|_{Q}$, $\forall x\in E$. In fact, let
$f\in E^{\ast}_{c}$, then there exists some $B_{Q}(\varepsilon)$
as in Proposition 1.5 such that
$f(B_{Q}(\varepsilon))\subset\{\xi\in L^{0}({\cal
F},K)\,|\,|\xi|\leq 1\}$, so that we can have
$|f(\frac{\varepsilon}{\|x\|_Q+1/n}x)|\leq 1,\forall x\in E$ and
$n\in N$, which means that $|f(x)|\leq\xi\|x\|_Q,\forall x\in E$,
where $\xi=1/\varepsilon$, and the converse is obvious. Lemma 4.1
of [25] is a special case of this result when $(E,{\cal P})$ is an
$RN$ module.

It is, however, not trivial to characterize an element in
$E^{\ast}_{\varepsilon,\lambda}$, as shown in [24]:

{\noindent\bf Proposition 1.6 ([24])}. A random linear functional $f:E\rightarrow L^{0}({\cal F},K)\in
E^{\ast}_{\varepsilon,\lambda}$ iff there are a countable partition
$\{A_{n}\,|\,n\in N\}$ of $\Omega$ to ${\cal F}$, a sequence
$\{\xi_{n}~|~n\in N\}$ in $L^{0}_+$ and a sequence $\{Q_{n}~|~n\in
N\}$ of finite subfamilies of $\cal P$ such that
$|f(x)|\leqslant\sum_{n\geq1}\tilde{I}_{A_{n}}\cdot\xi_{n}\cdot\|x\|_{Q_{n}}$,
$\forall x\in E$, in which case if, let
$\xi=\sum_{n\geq1}\tilde{I}_{A_{n}}\cdot\xi_{n}$ and
$\|x\|=\sum_{n\geq1}\tilde{I}_{A_{n}}\cdot\|x\|_{Q_{n}}$, $\forall
x\in E$, then $|f(x)|\leq\xi\cdot\|x\|,\forall x\in E$.

Proposition 1.6 shows that $E^{\ast}_{\varepsilon,\lambda}\supset
E^{\ast}_{c}$, and $E^{\ast}_{\varepsilon,\lambda}= E^{\ast}_{c}$ if
$\mathcal {P}$ has the countable concatenation property, in
particular $E^{\ast}_{\varepsilon,\lambda}= E^{\ast}_{c}$ for any
$RN$ module $(E,\|\cdot\|)$.

$E^{\ast}_{\varepsilon,\lambda}$ and $E^{\ast}_{c}$ are both left modules
over the algebra $L^{0}({\cal F},K)$ by $(\xi\cdot
f)(x)=\xi\cdot(f(x))$, $\forall\xi\in L^{0}({\cal F},K),f\in
E^{\ast}_{\varepsilon,\lambda}$ or $E^{\ast}_{c}$, and $x\in E$, the
following Definition 1.7 shows that
$E^{\ast}_{\varepsilon,\lambda}=E^{\ast}_{c}$ is still a complete
$RN$ module for an $RN$ module $(E,\|\cdot\|)$, see Section 1.4 for
completeness.

{\noindent\bf Definition 1.7 ([11])}. Let $(E,\|\cdot\|)$ be an $RN$
module over $K$ with base $(\Omega,\mathscr{A},\mu)$ and
$E^{\ast}:=E^{\ast}_{\varepsilon,\lambda}=E^{\ast}_{c}$. Define
$\|\cdot\|^{\ast}:E^{\ast}\rightarrow L^{0}_+$ by
$\|f\|^{\ast}=\wedge\{\xi\in
L^{0}_+\,|\,|f(x)|\leq\xi\cdot\|x\|,\forall x\in E\}$, then
$(E^{\ast},\|\cdot\|^{\ast})$ is an $RN$ module over $K$ with base
$(\Omega,{\cal F},P)$, called the random conjugate space of $(E,
\|\cdot\|)$, which is complete under both ${\cal
T}_{\varepsilon,\lambda}$ and ${\cal T}_{c}$ (notice:
$\|f\|^{\ast}=\vee\{|f(x)|\,|\,x\in E$ and $\|x\|\leq 1\}$, cf.
[19]).

\subsection*{1.3. The countable concatenation closure of a set and hyperplane separation
theorems.}
\quad\, The hyperplane separation Theorem 2.6 of [4] is
peculiar to ${\cal T}_{c}$ since the gauge function perfectly
matches ${\cal T}_{c}$, it is impossible to present it under
${\cal T}_{\varepsilon, \lambda}$, as said in [22]. The hyperplane
separation Theorem 2.8 of [4] and Theorem 3.1 of [22] are both a
random generalization of the famous Mazur's theorem, it is not
difficult to prove that our Theorem 3.1 of [22] implies the
separation Theorem 2.8 of [4], in particular, our Theorem 3.1 of
[22] allows a kind of separation with an arbitrary probability,
and thus it can, like the classical Mazur's theorem, implies that
an $L^{0}-$convex subset(an $M$-convex subset in terms of [11, 19,
22]) is ${\cal T}_{\varepsilon,\lambda}-$closed iff it is random
weakly closed (see Corollary 3.4 of [22]). But, the separation
Theorem 2.8 of [4] can not derive such a kind of result, since it
was only given in a form of separation with probability one.
Though the hyperplane separation Theorems 2.6 and 2.8 of [4] have
succeeded in convex analysis for locally $L^{0}-$convex modules,
we would like to generalize them to a more general form allowing a
kind of separation with an arbitrary probability. In the process,
the notion of a countable concatenation closure of a set will play
a key role, see Section 3 of this paper.
\subsection*{ 1.4. Completeness}
\quad\, ${\cal T}_{\varepsilon,\lambda}-$completeness and ${\cal
T}_{c}-$completeness are very important both in convex analysis
for locally $L^{0}-$convex modules, cf.[25], and in the deep study
of random conjugate spaces, cf.[16 ,18 ,19 ,23]. It is easy to
observe that ${\cal T}_{\varepsilon,\lambda}-$completeness of a
random locally convex module $(E,\cal P)$ implies its ${\cal
T}_{c}-$completeness since ${\cal T}_{c}$ has a local base ${\cal
U}_{\theta}=\{B_{Q}(\varepsilon)~|~Q\subset\cal P$ finite and
$\varepsilon\in L_{++}^{0}\}$ as in Proposition 1.5 such that each
$B_{Q}(\varepsilon)$ is ${\cal T}_{\varepsilon, \lambda}-$closed,
it follows immediately from this observation that $E^{\ast}$ and
$L^{0}({\cal F},K)$ are ${\cal T}_{c}-$complete for an $RN$ module
$(E,\|\cdot\|)$ since it is very easy to verify that they are
${\cal T}_{\varepsilon,\lambda}-$complete, cf.[9, 18]. But it is a
delicate matter whether ${\cal T}_{c}-$completeness also implies
${\cal T}_{\varepsilon,\lambda}-$completeness, since one will find
that it is not easy to prove this. In the final part of Section 3,
we combine the very interesting Theorem 3.18 of this paper and the
usual completion skill so that the following statement can be
obtained: ${\cal T}_{c}-$completeness of a random locally convex
module with the countable concatenation property implies ${\cal
T}_{\varepsilon,\lambda}-$completeness, where the countable
concatenation property is different from the one in the sense of
[4], and it is easy to verify that all the random normed modules
currently used in conditional risk measures, cf.[25], have the new
property.

\subsection*{1.5. The theory of random conjugate spaces of random normed
modules under ${\cal T}_{c}$} \quad\, In the past ten years, the
theory of random conjugate spaces of random normed modules under
${\cal T}_{\varepsilon,\lambda}$ has been quite deep and systematic,
cf.[8, 10, 16, 19, 23], and it is comparatively difficult to
establish the results in [8, 10, 16, 19, 23]. But the corresponding
theory under ${\cal T}_{c}$ is still at the beginning stage, e.g.,
the representation theorems of random conjugate spaces of the only
$L^{0}({\cal F},K^{d})$ and $L^{p}_{{\cal F}}({\cal E})$ were given
in [25]. Based on the above discussion of completeness, in Section 4
we prove that all the results obtained in [8 ,10 ,19 ,23] are still
valid under ${\cal T}_{c}$, but those in [16] is no longer valid
under ${\cal T}_{c}$ since ${\cal T}_{c}$ is too strong. Further,
motivated by the idea of Filipovi\'{c}, kupper and vogelpoth's
constructing $L^{p}_{{\cal F}}({\cal E})$ and representing its
random conjugate space, we construct $L^{p}_{{\cal F}}(E)$ and
represent $(L^{p}_{{\cal F}}(E))^{\ast}$ as $L^{q}_{{\cal
F}}(E^{\ast})$ in an isomorphically isometric manner, where ${\cal
F}\subset\cal E,{\cal F}$ and $\cal E$ are both the
$\sigma-$algebras over $\Omega$, and $E$ is an arbitrary random
normed module over $K$ with base $(\Omega,\cal E,\cal P)$.
Specially, take $E=L^{0}({\cal E},R)$, then we have $L^{p}_{{\cal
F}}({E})=L^{p}_{{\cal F}}({\cal E})$, and thus generalize and
strengthen the corresponding representation theorem of [4].
Consequently the whole Section 4 has established, on a large scale,
the theory of random conjugate spaces of random normed modules under
${\cal T}_{c}$, which considerably
enriches financial applications of random normed modules.

By the way, although there have been several proofs of the known
Hahn-Banach extension theorem for $L^0-$linear functions, cf.[2, 4,
30], each of them has to face difficulties in the existence of ``one
step extension", and thus comparatively complicated. In Section 2,
we first give an extremely simple proof, which avoids the above
difficulties by reducing an extension problem for an $L^0-$linear
function on an $L^0-$submodule to one for a random linear functional
on a linear subspace. We also consider the Hahn-Banach extension
theorem for a random linear functional on a complex subspace as well
as an $L^0(\mathcal{F}, C)-$linear function. Finally, we conclude
Section 2 with the Hah-Banach extension theorems for continuous
$L^0-$linear functions.

\section*{2. Hahn-Banach extension theorems}

\subsection*{2.1. Hahn-Banach extension theorems for random linear
functionals}

\quad\, Let us first recall from [13]: let $X$ be a linear space
over $K$, then a linear operator $f$ from $X$ to $L^0(\mathcal{F},
K)$ is called \emph{a random linear functional}  on $X$; A mapping
 $p:X\rightarrow L^{0}_+$ is called \emph{a random seminorm } on $X$ if
 it satisfies the following :

 $(1) \  \  \ p(\alpha x)=|\alpha| p(x), \forall \alpha \in K $ and $ x\in
 X;$

 $(2) \ \ \ p(x+y)\leq p(x)+p(y),   \forall x,  y\in X.
$

Let $X$ be a real linear space. A mapping $p:X\rightarrow
L^{0}(\mathcal{F}, R)$ is called \emph{a random sublinear functional
} if it satisfies the above (2) and the following:

 $(3)\ \ \ p(\alpha x)=\alpha p(x), \forall \alpha \geq 0\   \hbox{and}\ x\in
 X;$

 The main results of the subsection are stated as follows:

 {\noindent \bf{Theorem 2.1.} } Let $E$ be a real linear space, $M\subset E$ a linear
 subspace, $f:M\rightarrow L^0(\mathcal{F}, R)$ a random
 linear functional and $p:E\rightarrow L^0(\mathcal{F}, R)$ a
 random sublinear functional such that $f(x)\leq p(x)$, $\forall x\in
 M$. Then there exists a random linear functional\ $g:E\rightarrow L^0(\mathcal{F}, R)$
such that $g$ extends $f$ and $g(x)\leq p(x)$,  \ $\forall x\in E$.

{\noindent \bf{Theorem 2.2.}} \ Let $E$ be a complex linear space,
$M\subset E$ a complex linear subspace, $f:M\rightarrow
L^0(\mathcal{F}, C)$  a random linear functional and
$p:E\rightarrow L^0_{+}$ a random seminorm such that $|f(x)|\leq
p(x)$, $\forall x\in M.$ Then there exists a random linear
functional\ $g:E\rightarrow L^0(\mathcal{F}, C)$ such that $g$
extends $f$ and $|g(x)|\leq p(x)$, \ $\forall x\in E$.

The proofs of Theorems 2.1 and 2.2 are given in [5, 6, 7, 13] but
in an indirect way, see the survey paper [13] for details. In fact,
Theorem 2.1 is known in [2, 30] since $R$ is, of course, an
ordered ring, and its proof is only a copy of the classical
Hahn-Banach extension theorem for a real linear functional by
replacing the order-completeness of $R$ with the one of
$L^0(\mathcal{F}, R)$. But the proof of Theorem 2.2 is somewhat
different from its classical prototype, since $E$ is not an
$L^0(\mathcal{F}, C)-$module. The idea of proof comes from [5], we
give a direct proof of Theorem 2.2 for the reader's convenience.

{\noindent\bf Proof of Theorem 2.2.} Let\ $f_{1}:M\rightarrow
L^0(\mathcal{F}, R)$ be the real part of $f$, then it is easy to see
that $f(x)=f_1(x)-if_1(ix)$, $\forall x \in M$. Then it is clear
that $f_1$ satisfies the following:
$$f_1(x)\leq p(x), \ \forall x\in M \eqno(2.1)$$

Regarding $M$ and $E$ as a real linear spaces, then Theorem 2.1
yields a random linear functional $g_{1}:E\rightarrow
L^0(\mathcal{F},R)$ such that $g_1$ extends $f_1$ and $g_1$
satisfies:
$$g_1(x)\leq p(x), \ \forall x\in E \eqno(2.2)$$

Define $g:E\rightarrow L^0(\mathcal{F},C)$ by
$g(x)=g_1(x)-ig_1(ix)$, $\forall x\in E$, then it is easy to check
that $g$ is a random linear functional extending $f$, we want to
prove that $g$ satisfies the following:
$$|g(x)|\leq p(x), \ \forall x\in E \eqno(2.3)$$

Given $x\in E$, let $g_1^{0}(x), \ g^{0}(x)$ and $p^{0}(x)$ be
arbitrarily chosen representatives of $g_1(x), \ g(x)$ and $p(x)$,
respectively. Then we have the following relations:
$$
g_1^{0}(x)(\omega) \leq p^{0}(x)(\omega), \ a.s., \ \forall x\in E,
\eqno(2.4)$$
$$g_1^{0}(x)(\omega)=Re(g^{0}(x)(\omega)), \ a.s.,\ \forall x\in E, \eqno(2.5)$$
$$g^{0}(x+y)(\omega)=g^{0}(x)(\omega)+g^{0}(y)(\omega), \ a.s., \ \forall x, \ y \in
E,\eqno(2.6)$$
 $$g^{0}(\alpha x)(\omega)=\alpha\cdot g^{0}(x)(\omega), \ a.s., \
\forall \alpha \in C \ and \  x\in E, \eqno(2.7)
$$
$$p^{0}(\alpha x)(\omega)=|\alpha|\cdot p^{0}(x)(\omega), \ a.s.,\
\forall \alpha \in C \ and \  x\in E. \eqno(2.8)
$$

It follows immediately from (2.7), (2.5), (2.4) and (2.8) that
$Re(\alpha\cdot g^{0}(x)(\omega))  = Re( g^{0}(\alpha x)(\omega))$
$=g_1^{0}(\alpha x)(\omega)\leq p^{0}(\alpha
x)(\omega)=|\alpha|\cdot p^{0}( x)(\omega), \ a.s., \ \forall
\alpha \in C \ and \  x\in E.$

Now, let $x$ be an arbitrary but fixed element of $E$ and
$\{c_{n}\,|\,n\in N\}$  a countable subset dense in $C$. Then there
exists $A_{n}(x)\in \mathcal{F}$ for each $n\in N$ such that
$P(A_{n}(x))=1$ and $Re(c_{n}\cdot g^{0}(x)(\omega))\leq
|c_{n}|\cdot p^{0}(x)(\omega), \ \forall \omega\in A_{n}(x)$.

Let $A(x)=\bigcap _{n\in N}A_{n}(x)$. Then  $A(x)\in \mathcal{F}, \
P(A(x))=1$ and the following is also true:
$$Re(\alpha \cdot g^{0}(x)(\omega))\leq |\alpha|\cdot p^{0}(x)(\omega), \ \forall
\omega\in A(x) \ and \ \alpha \in C, \eqno(2.9)
$$
since $\{c_{n}\,|\,n\in N\}$ is dense in $C$.

Given $\omega\in \Omega$, let $\theta (\omega)$ be the principal
argument of $g^{0}(x)(\omega)$, then $\theta$ is a random variable
and $|g^{0}(x)|(\omega)=e^{-i\theta (\omega)}\cdot g^{0}(x)(\omega),
\ \forall \omega\in \Omega.$

It follows immediately from (2.9) that
$|g^{0}(x)|(\omega)=Re(|g^{0}(x)|(\omega))=Re(e^{-i\theta
(\omega)}\cdot g^{0}(x)(\omega))\leq |e^{-i\theta (\omega)}|\cdot
p^{0}(x)(\omega)= p^{0}(x)(\omega), \ \forall \omega\in A(x)$. Thus,
we have $|g^{0}(x)|(\omega)\leq p^{0}(x)(\omega),a.s., $ since
$P(A(x))=1$. Namely, $|g(x)|\leq p(x),$ which comes from the
definition of the ordering $\leq$.

Finally, since $x$ is arbitrary, then (2.3) has been proved.
$\square$

{\noindent\bf{Remark 2.3. }}Let $\mathscr{L}^0(\mathcal{F}, K)$ be
the linear space of $K-$valued $ \mathcal{F}-$measurable random
variables on $(\Omega, \ \mathcal{F}, \ P)$ under the ordinary
pointwise scalar multiplication and addition operations, it is
easy to see that $L^{0}(\mathcal{F}, K)$ is just the quotient
space of $\mathscr{L}^0(\mathcal{F}, K)$ under the equivalence
relation $\sim : \xi \sim \eta$ iff $\xi (\omega)=\eta(\omega), \
a.s.$ Namely, when elements equal a.s. are identified,
$\mathscr{L}^0(\mathcal{F}, K)$ is exactly $L^0(\mathcal{F}, K)$.
But some situations do not allow us to regard elements equal a.s.
in $\mathscr{L}^0(\mathcal{F}, K)$ as identified: for example, let
$g^{0}:E\rightarrow \mathscr{L}^0(\mathcal{F}, K)$ (by replacing
$C$ with $K$) be the mapping as defined in the proof of Theorem
2.2. Conditions (2.6) and (2.7) shows that $g^{0}$ is a kind of
random linear operator (precisely, a random linear functional) as
defined in [27]. If for each $\omega\in\Omega$,
$g^{0}(\cdot)(\omega):E\rightarrow K$ is a linear functional, then
$g^{0}$ is called sample-linear, in which case $g^{0}:
E\rightarrow \mathscr{L}^0(\mathcal{F}, K)$ is an ordinary linear
operator! To obtain a sample-linear random functional from a
random linear functional, we often employ the theory of lifting
property, cf. [28, 29], in all the cases when we consider a
sample$-$linear random functional, we can not, of course, identify
$L^{0}({\cal F},K)$ with $\mathscr {L}^{0}({\cal F},K)$.
Sometimes, the domain of $g^{0}$ is not the whole $E$ but a random
subset in $E$, we are forced to use measurable selection theorems
of [31], cf.[18]. Since the theories of [27, 28, 31] all require
the completeness of $\mathcal{F}$ with respect to $P$, but
$\hat{\mathcal{F}}=$ the completion of $\mathcal {F}$ with respect
to $P$ is, of course, complete, consequently we used to employ
$L(P, K)$ in [8-24], namely $L^0(\hat{\mathcal {F}}, K)$, so
that we need not assume $\mathcal {F}$ to be complete.

{\noindent\bf{Remark 2.4.}} Before 1996, we employed the original
definition of a random normed space, and thus the Hahn-Banach
extension theorem for a random linear functional was given for such
as $g^{0,}$s in Remark 2.3, in [5, 6, 7], which leads to a kind of
theory of random conjugate space of a random normed space, see[13,
Section 3]. After 1996, in particular after 1999, the author gave a
new version of a random normed space in [11], we began to consider
the Hahn-Banach extension theorem for a random linear functional $g$
as in Theorems 2.1 and 2.2, which leads to the current frequently
used theory of a random conjugate space, cf. [13, Section 4].

\subsection*{2.2. Hahn-Banach extension theorems for $L^0-$linear
functions}

\quad\, Let $E$ be a left module over the algebra $L^0(\mathcal
{F}, K)$ (equivalently, the ring $L^0(\mathcal {F}, K)$, see
Remark 1.2 for details), a module homomorphism $f:E\rightarrow
L^0(\mathcal {F}, K)$, is called an $L^0$(or $L^0(\mathcal {F},
K))-$linear function. If $K=R$, then a mapping $p:E\rightarrow
L^0(\mathcal {F}, R)$ is called an $L^0-$sublinear function if it
satisfies the following:

\hspace{2mm}(\lowercase\expandafter{\romannumeral 1}) $p(\xi x)=\xi
p(x)$, $\forall\xi\in L_{+}^{0}$ and $x\in E$;

(\lowercase\expandafter{\romannumeral 2}) $p(x+y)\leq p(x)+p(y)$,
$\forall x,y\in E$.

The main results of this subsection are known, see[2, 4, 30], but we
give them extremely simple proofs. They are stated as follows:

{\noindent\bf{Theorem 2.5.}} Let $E$ be a left module over the
algebra $L^0(\mathcal {F}, R)$, $M\subset E$ an $L^0(\mathcal {F},
R)-$submodule, $f:M\rightarrow L^0(\mathcal {F}, R)$ an $L^0-$linear
function and $p:E\rightarrow L^0(\mathcal {F}, R)$ an
$L^0-$sublinear function such that $f(x)\leq p(x)$, $\forall x\in
M$.Then there exists an $L^{0}-$linear function $g:E\rightarrow
L^0(\mathcal {F}, R)$ such that $g$ extends $f$ and $g(x)\leq p(x),
\ \forall x\in E$.

{\noindent\bf{Theorem 2.6.}}  Let $E$ be a left module over the
algebra $L^0(\mathcal {F}, C)$, $M\subset E$ an $L^0(\mathcal {F},
C)-$submodule, $f:M\rightarrow L^0(\mathcal {F}, C)$ an $L^0-$linear
function and $p:E\rightarrow L^0_{+}$ an $L^0-$seminorm such that
$|f(x)|\leq p(x)$, $\forall x\in M$. Then there exists an
$L^{0}-$linear function $g:E\rightarrow L^0(\mathcal {F}, C)$ such
that $g$ extends $f$ and $|g(x)|\leq p(x), \ \forall x\in E$.

Theorem 2.6 follows immediately from Theorem 2.5, since $E$ is a
module and the method used to prove its classical prototype will
still be feasible. The following simple lemma can lead to an
extremely simple proof of Theorem 2.5.

{\noindent\bf{Lemma 2.7.}} Let $E$ be a left module over the algebra
$L^0(\mathcal {F}, R)$, $f:E\rightarrow L^0(\mathcal {F}, R)$ a
random linear functional and $p:E\rightarrow L^0(\mathcal {F}, R)$
an $L^0-$sublinear function such that $f(x)\leq p(x)$, $\forall x\in
E$. Then $f$ is an $L^{0}-$linear function. If $R$ is replaced by
$C$ and $p$ is an $L^0-$seminorm such that $|f(x)|\leq p(x)$,
$\forall x\in E$, then $f$ is also an $L^{0}-$linear function.

{\noindent \bf{Proof.}} We only give the proof of the first part,
since the second one is similar.

Since $f$ is linear, it suffices to prove that $f(\xi x)=\xi f(x), \
\forall \xi \in L^{0}_{+}$ and $x\in E$.

Let $x\in E$ be fixed. For any $\xi \in L^{0}_{+}$, since there
exists a sequence $\{\xi_{n}\,|\,n\in N \}$ of simple elements in
$L^{0}_{+}$ such that $\{\xi_{n}\,|\,n\in N \}$ converges a.s. to
$\xi$ in a nondecreasing way, and since $|f(\xi x)-f(\xi_{n}
x)|=|f[(\xi -\xi_{n})x]|\leq (\xi -\xi_{n})(|p(x)|\vee |p(-x)|)$, it
also suffices to prove that $f({\tilde{I_{A}}} x)=
{\tilde{I_{A}}}\cdot f(x), \
   \forall A\in \mathcal {F}$.

Since $-p(-\tilde{I_{A}} x)\leq f(\tilde{I_{A}} x)\leq
p(\tilde{I_{A}}x)$, namely $-\tilde{I_{A}}p(- x)\leq f(\tilde{I_{A}}
x)\leq \tilde{I_{A}} p(x)$, we have
$\tilde{I_{A^{c}}}f(\tilde{I_{A}} x)=0, \ \forall A\in \mathcal
{F}$, where $A^{c}=\Omega\setminus A$. Then we have that
$f(\tilde{I_{A}} x)=\tilde{I_{A}}f(\tilde{I_{A}}
x)+\tilde{I_{A^{c}}}f(\tilde{I_{A}} x)=\tilde{I_{A}}f(\tilde{I_{A}}
x)=\tilde{I_{A}}f(\tilde{I_{A}} x)+\tilde{I_{A}}f(\tilde{I_{A^{c}}}
x)=\tilde{I_{A}}f(\tilde{I_{A}}
x+\tilde{I_{A^{c}}}x)=\tilde{I_{A}}f(x)$. $\square$

We can now prove Theorem 2.5 in an extremely simple way.

{\noindent\bf Proof of Theorem 2.5.} Applying Theorem 2.1 to $f$ and
$p$ yields a random linear functional $g:E\rightarrow L^{0}(\mathcal
{F},R)$ such that $g$ extends $f$ and $f(x)\leq p(x), \ \forall x\in
E$. Since $E$ is a left module over the algebra $L^{0}(\mathcal
{F},R)$ and $p$ is an $L^{0}-$sublinear function, Lemma 2.7 shows
that $g$ is again an $L^{0}-$linear function. $\square$

{\noindent \bf{Remark 2.8.}} The idea of proof of Theorem 2.5 comes
from [7], which is also used in [20, 24].

\subsection*{2.3. Hahn-Banach extension theorems for continuous
$L^{0}-$linear functions in random locally convex modules under the
two kinds of topologies.}

{\noindent \bf{Theorem 2.9.}} Let $(E,\mathcal{P})$ be a random
locally convex module over $K$ with base $(\Omega,\mathcal{F},P)$
and $M$ an
 $L^{0}(\mathcal{F},K)-$submodule of $E$. Then we have the following:

(1)\, Every continuous $L^{0}-$linear function $f$ from $(M,
\mathcal {T}_{c})$ to $(L^{0}(\mathcal{F},K), \mathcal {T}_{c})$
admits a continuous $L^{0}-$linear extension $g$ from $(E, \mathcal
{T}_{c})$ to $(L^{0}(\mathcal{F},K), \mathcal {T}_{c})$;

(2)\, Every continuous $L^{0}-$linear function $f$ from $(M,
\mathcal {T}_{\varepsilon, \lambda})$ to $(L^0(\mathcal {F}, K),
\mathcal {T}_{\varepsilon, \lambda})$ admits a continuous
$L^0-$linear extension $g$ from $(E, \mathcal {T}_{\varepsilon,
\lambda})$ to $(L^0(\mathcal {F}, K), \mathcal {T}_{\varepsilon,
\lambda})$;

(3)\, If $(E,\mathcal {P})$ is an $RN$ module $(E, \|\cdot\|)$, then
$g$ in both (1) and (2) can be required to be such that
$\|g\|^*=\|f\|^*$, namely $g$ is a random-norm preserving extension.

{\noindent \bf{ Proof.}} (1). Let $\mathcal {P}_M=\{\|\cdot\|_M\mid
\|\cdot\|\in\mathcal {P}\}$, where $\|\cdot\|_M$ is the restriction
of $\|\cdot\|$ to $M$, then $(M, \mathcal {P}_M)$ is still a random
locally convex module. Since $f$ is a continuous $L^0-$linear
function from $(M, \mathcal {T}_c)$ to $(L^0(\mathcal {F}, K),
\mathcal {T}_{c})$, there exists a finite subfamily $Q$ of $\mathcal
{P}$ and $\xi\in L^0_+$ such that $|f(x)|\leq \xi\|x\|_Q$, $\forall
x\in M$. Then Theorem 2.5 and Theorem
2.6 jointly yield an $L^0-$linear function $g:E\rightarrow
L^0(\mathcal {F}, K)$ such that $g$ extends $f$ and $|g(x)|\leq
\xi\|x\|_Q$, $\forall x\in E$. Of course, $g$ satisfies the requirement of (1).

(2). It follows immediately from Proposition 1.6 that there exists
$\xi\in L^0_+$, a countable partition $\{A_n\mid n\in N\}$ of
$\Omega$ to $\mathcal {F}$ and a sequence $\{Q_n\mid n\in N\}$ of
finite subfamilies of $\mathcal {P}$ such that $|f(x)|\leq
\xi\|x\|$, $\forall x\in M$, where $\|\cdot\|:E\rightarrow L^0_+$
is given by $\|x\|=\sum_{n\geq 1}\widetilde{I}_{A_n}\|x\|_{Q_n}$,
$\forall x\in E$.

Obviously, $\|\cdot\|$ is an $L^0-$seminorm, then Theorems 2.5 and
2.6 jointly yield an $L^0-$linear function $g:E\rightarrow
L^0(\mathcal {F}, K)$ such that $|g(x)|\leq \xi\|x\|$, $\forall
x\in E$. Further, since $\sum_{n\geq 1}P(A_n)=P(\Omega)=1$, then
$P(A_n)\rightarrow 0$ as $n\rightarrow\infty$, and hence
$\|\cdot\|$ is continuous from $(E, \mathcal {T}_{\varepsilon,
\lambda})$ to $(L^0_+, \mathcal {T}_{\varepsilon, \lambda})$,
which implies that $g$ is also continuous.

(3). Let $p(x)=\|f\|^*\|x\|$, $\forall x\in E$, then Theorems 2.5
and 2.6 again jointly complete the proof. \quad $\square$

{\noindent \bf{Remark 2.10.}} Theorem 2.9 is given to tidy up the
three diverse results in it. (2) of Theorem 2.9 is already given in
[24] as a corollary of Proposition 1.6 but the proof in [24] did not
use Theorems 2.5 and 2.6 but used the idea of proof of Lemma 2.7;
(3) of Theorem 2.9 has been known for at least 20 years, as the
Hahn-Banach extension theorem for a.s. bounded random linear
functionals defined on a random normed space, which was given in
[5]; Since such a $L^0-$seminorm $\|\cdot\|$ as in the proof of (2)
of Theorem 2.9 is always $\mathcal {T}_{\varepsilon,
\lambda}-$continuous, we can, without loss of generality, assume
that $\mathcal {P}$ has the countable concatenation property if we
only consider $\mathcal {T}_{\varepsilon, \lambda}$.

\section*{3. Hyperplane separation theorems }

\subsection*{ 3.1. The countable concatenation hull, the
countable concatenation property and relations between hyperplane
separation theorems currently available.}

 \quad\, Filipovi\'{c}, Kupper and Vogelpoth introduced in [4] the two
countable concatenation properties, one is relative to topology and
the other is relative to the family of $L^0-$seminorms, in
particular the two are essentially the same one for a random locally
convex module $(E, \mathcal {P})$ endowed with $\mathcal {T}_c$.
Thus neither of the two has anything to do with the $L^0-$module $E$
itself. The main results in Sections 3 and 4 show that there is
another kind of countable concatenation property, which is concerned
with the $L^0-$module $E$ itself and is very important for the
theory of locally $L^0-$convex module. To introduce it, we give the
following:

First, we make the following convention that all the $L^0(\mathcal
{F}, K)-$modules $E$ in the sequel of this paper satisfy the
property: For any two elements $x$ and $y\in E$, if there exists a
countable partition $\{A_n\mid n\in N\}$ of $\Omega$ to $\mathcal
{F}$ such that $\widetilde{I}_{A_n}x=\widetilde{I}_{A_n}y$ for each
$n\in N$, then $x=y$, where $\widetilde{I}_{A}x$ is called the
$A-$stratification of $x$ for any $A\in\mathcal {F}$ (see [16, 22,
24] for some discussions of stratification structure). Clearly, a
random locally convex module $(E, \mathcal {P})$ always satisfies
the convention, since $\vee\{\|x-y\|\mid\|\cdot\|\in {\cal P}\}=\sum_{n\geq
1}\widetilde{I}_{A_n}(\vee\{\|x-y\|\mid\|\cdot\|\in {\cal P}\})=\sum_{n\geq
1}(\vee\{\|\widetilde{I}_{A_n}x-\widetilde{I}_{A_n}y\|\mid\|\cdot\|\in
{\cal P}\})=0$ if $\widetilde{I}_{A_n}x=\widetilde{I}_{A_n}y, \forall n\in
N$, then we have $x=y$ by the definition of a random locally convex
module!

{\noindent\bf Definition 3.1.} Let $E$ be a left module over the
algebra $L^0(\mathcal {F}, K)$. Such a formal sum $\sum_{n\geq
1}\widetilde{I}_{A_n}x_n$ for some countable partition $\{A_n\mid
n\in N\}$ of $\Omega$ to $\mathcal {F}$ and some sequence $\{x_n\mid
n\in N\}$ in $E$, is called a \emph{{countable concatenation}} of
$\{x_n\mid n\in N\}$ with respect to $\{A_n\mid n\in N\}$.
Furthermore a countable concatenation $\sum_{n\geq
1}\widetilde{I}_{A_n}x_n$ is \emph{{well defined or $\sum_{n\geq
1}\widetilde{I}_{A_n}x_n\in E$}} if there is $x\in E$ such that
$\widetilde{I}_{A_n}x=\widetilde{I}_{A_n}x_n, \forall n\in N$
(Clearly, $x$ is unique, in which case we write $x=\sum_{n\geq
1}\widetilde{I}_{A_n}x_n$). A subset $G$ of $E$ is called
 \emph {having the countable concatenation property} if
every countable
concatenation $\sum_{n\geq 1}\tilde{I}_{A_n}x_n$ with $x_n\in
 G$ for each $n\in N$ still belongs to $G$, namely $\sum_{n\geq
 1}\tilde{I}_{A_n}x_n$ is well defined and there exists $x\in G$
 such that $x=\sum_{n\geq 1}\tilde{I}_{A_n}x_n$.

{\noindent\bf Definition 3.2.}  \  Let $E$ be a left module over the
algebra $L^{0}({\cal F},K)$ and $\{G_n\,|\,n\in N\}$ a sequence of
subsets of $E$. The set of well defined countable concatenations
$\sum_{n\geq 1}\tilde{I}_{A_n}x_n$ with $x_n\in
 G_n$ for each $n\in N$ is called {\it the countable concatenation
 hull} of the sequence $\{G_n\,|\,n\in N\}$, denoted by $H_{cc}(\{G_n\,|\,n\in
 N\})$. In particular when $G_n=G$ for each $n\in N$, $H_{cc}(G)$,
 denoting $H_{cc}(\{G_n\,|\,n\in
 N\})$,
 is called {\it the countable concatenation
 hull} of the subset $G$ (Clearly $H_{cc}(G)\supset G,\forall\, G\subset E$, and
 $H_{cc}(E)=E$).

 Definition 3.1 leads to a notion of $E$ having the countable
 concatenation property as a subset of itself. To remove any
 possible confusions, we will reserve the terminologies of [4] in
 the following equivalent way:

{\noindent\bf Definition 3.3([4]).}  Let $(E,{\cal T})$ be a
topological module over the topological ring $(L^{0}({\cal
F},K),{\cal T}_c)$. ${\cal T}$ is called {\it having the countable
 concatenation property } if $H_{cc}(\{U_n\,|\,n\in
 N\})$ is again a neighborhood of $\theta$ (the null element of $E$)
 for every sequence $\{U_n\,|\,n\in
 N\}$ of neighborhoods of $\theta$. Let $(E,{\cal P})$ be a random
 locally convex module over $K$ with base $(\Omega,{\cal F},P)$,
 ${\cal P}$ is called {\it having the countable
 concatenation property} if $\sum_{n\geq 1}\tilde{I}_{A_n}\cdot
 \|\cdot\|_{Q_n}$ still belongs to ${\cal P}$ for any countable
 partition $\{A_n\,|\,n\in
 N\}$ of $\Omega$ to ${\cal F}$ and any sequence $\{Q_n\,|\,n\in
 N\}$ of finite subfamilies of ${\cal P}$.

 Example 3.4 below, adopted from [4], exhibits that it is very
 necessary to introduce the notion of a well defined countable
 concatenation.

 {\noindent\bf Example 3.4.} For an $L^{0}({\cal F},K)-$module
 $E$, $M\subset E$ a subset of $E$, the
 $L^{0}({\cal F},K)-$submodule generated by $M$ is denoted by
 $Span_{L^{0}}(M)$, namely $Span_{L^{0}}(M)=\{\sum_{i=1}^{n}\xi_{i}x_{i}\,|\,\xi_{i}\in L^{0}({\cal F},K), x_i$ $\in M, 1\leq i\leq n$ and $n\in
 N\}$. Take $E=L^{0}({\cal F},R)$, $M=\{\tilde{I}_{[1-2^{-(n-1)},1-2^{-n}]}\,|\,n\in
 N\}$ and $F=Span_{L^{0}}(M)$, then it is easy to see that $\sum_{n=1}^{\infty}2^{n}\tilde{I}_{[1-2^{-(n-1)},1-2^{-n}]}\in E\backslash
 F$ and the sequence $\{\sum_{n=1}^{k}2^{n}\tilde{I}_{[1-2^{-(n-1)},1-2^{-n}]}\,|\,k\in
 N\}$ in $F$ is not a ${\cal T}_c-$Cauchy sequence but a
 ${\cal T}_{\varepsilon,\lambda}-$Cauchy sequence which has no limit
 in $(F,{\cal T}_{\varepsilon,\lambda})$. Thus the countable
 concatenation
 $\sum_{n=1}^{\infty}2^{n}\tilde{I}_{[1-2^{-(n-1)},1-2^{-n}]}$
 is not well defined in both $(F,{\cal T}_c)$ and
 $(F,{\cal T}_{\varepsilon,\lambda})$ in any ways, namely neither of
$(F,{\cal T}_c)$ and $(F,{\cal T}_{\varepsilon,\lambda})$ has the
countable
 concatenation property.

 However, the following examples show that all the
 $L^{0}({\cal F},K)-$modules important for
 financial applications have the countable
 concatenation property.

 {\noindent\bf Example 3.5.} $L^{p}_{{\cal F}}({\cal E})$ in
 [4,25] and the Orlicz type of $RN$ module $L^{\varphi}_{{\cal F}}({\cal
 E})$ in [25] both have the countable
 concatenation property, which is easily seen by the smooth property
 of the conditional expectation of random variables.

 {\noindent\bf Example 3.6.} Let $E$ be an
 $L^{0}({\cal F},K)-$module and $E^{\#}=\{f:E\rightarrow L^{0}({\cal F},K) \ | \ f$ is a module
 homomorphism$\}$, then it is easy to check that $E^{\#}$ has the countable
 concatenation property. Further, let $(E,{\cal P})$ be a random
 locally convex module over $K$ with base $(\Omega,{\cal F},P)$,
 $E_{c}^{\ast}$ and $E_{\varepsilon,\lambda}^{\ast}$ the corresponding
 random conjugate spaces, see Section 1.2. Then
 $E_{\varepsilon,\lambda}^{\ast}$ always has the countable
 concatenation property by Proposition 1.6, and $E_{c}^{\ast}$ has
 the countable concatenation property if ${\cal P}$ has the countable
 concatenation property by the paragraph in front of Proposition
 1.6.

 Example 3.7 below, in principle, surveys an extremely useful skill in the
 history of the development of a
 ${\cal T}_{\varepsilon,\lambda}-$complete random locally convex
 module, which played a crucial role in [8,10,16,17,19,22,23]. The
 central idea of this method is described as follows: To seek one
 desired element $x$ in a ${\cal T}_{\varepsilon,\lambda}-$complete random locally convex
 module $E$, we are forced to first find out a sequence $\{x_n\,|\,n\in
 N\}$ in $E$ such that each $x_n$ is the $A_n-$stratification of $x$
 and $\{A_n\,|\,n\in N\}$ exactly forms a countable partition of
 $\Omega$ to ${\cal F}$, since $\{\sum_{n=1}^{k}\tilde{I}_{A_n}x_n\,|\,k\in
 N\}$ is easily verified to be a
 ${\cal T}_{\varepsilon,\lambda}-$Cauchy sequence since $P(A_n)\rightarrow
 0$ as $n\rightarrow\infty$, whose limit is just the desired $x$!

{\noindent\bf Example 3.7.} Every ${\cal
T}_{\varepsilon,\lambda}-$complete random locally convex module
$(E,{\cal P})$ has the countable concatenation property: since for
every countable concatenation
$\sum_{n=1}^{\infty}\tilde{I}_{A_n}x_n, \
\{\sum_{n=1}^{k}\tilde{I}_{A_n}x_n\,|\,k\in N\}$ is a ${\cal
T}_{\varepsilon,\lambda}-$Cauchy sequence, and hence convergent to
some $x\in E$ so that $\tilde{I}_{A_n}x_n=\tilde{I}_{A_n}x,\forall
x\in E$, namely $\sum_{n\geq 1}\tilde{I}_{A_n}x_n$ is well defined.

Now, we will use the notions introduced above to study the relations
between hyperplane separation theorems currently available. To do
this, we first state the three hyperplane separation theorems, the
first two in [4] and the third in [22], in the following equivalent
ways.

{\noindent\bf Theorem 3.8([4, Theorem 2.6, Hyperplane separation
\uppercase\expandafter{\romannumeral 1}]). } Let $(E,{\cal P})$ be a
random locally convex module over $R$ with base $(\Omega,{\cal
F},P)$, $M$ and $G$ be $L^{0}-$convex subsets of $E$ such that $G$
is a nonempty ${\cal T}_c-$open subset. If $G$ and $M$ satisfy the
following:

$\tilde{I}_AM\,\cap\, \tilde{I}_AG=\emptyset$ for all $A\in{\cal F}$
with $P(A)>0,$ $\hfill{}$ (3.10)\\
then there exists an $f\in E_{c}^{\ast}$ such that

$f(x)<f(y)$ on $\Omega$ for all $x\in G$ and $y\in M.$$\hfill{}$ (3.11)\\

{\noindent\bf Remark 3.9.} If $R$ is replaced by $C$, then Theorem
3.8 still holds in the way: $(Ref)(x)<(Ref)(y)$ on $\Omega$ for all
$x\in G$ and $y\in M$, where $Ref:E\rightarrow L^{0}({\cal F},R)$
defined by $(Ref)(x)=Re(f(x)),\forall x\in E$, is the real part of
$f$, it is easy to see that $f(x)=(Ref)(x)-i(Ref)(ix),\forall x\in
E$.

{\noindent\bf Theorem 3.10([4,Theorem 2.8, Hyperplane separation
\uppercase\expandafter{\romannumeral 2}]).}\,   Let $(E,{\cal P})$
be a random locally convex module over $R$ with base $(\Omega,{\cal
F},P)$ such that ${\cal P}$ has the countable concatenation
property, $x\in E$ and $G\subset E$ a nonempty ${\cal T}_c-$closed
$L^{0}-$convex subset with the countable concatenation property. If
$x$ and $G$ satisfy the following:

$\tilde{I}_A\{x\}\,\cap\,\tilde{I}_AG=\emptyset$ for all
$A\in{\cal F}$ with $P(A)>0$,$\hfill{}$ (3.12)\\
then there exist an $f\in E_{c}^{\ast}$ and $\varepsilon\in
L^{0}_{++}$ such that

$f(x)>f(y)+\varepsilon$ on $\Omega$ for all $y\in G$.$\hfill{}$ (3.13)\\
Furthermore, if $R$ is replaced by $C$, then Theorem 3.10 is still
true in the way: $(Ref)(x)>(Ref)(y)+\varepsilon$ on $\Omega$ for all
$y\in G$.

{\noindent\bf Remark 3.11.} In [4], the original Theorem 2.8 of [4]
did not require $G$ to have the countable concatenation property,
but Lemma 3.17 below plays a key role in the proof of Theorem 2.8 of
[4] and the proof of Lemma 3.17 seems to ask $G$ to have such a
property. On the other hand, requiring $G$ to have the countable
concatenation property would not reduce the values of Theorem 3.10
and Lemma 3.17. Take, the important application of Theorem 3.10 to
Theorem 3.8 of [4] and the important application of Lemma 3.17 to
Lemma 3.10 of [4], for example, once the whole spaces $E$ in both
Lemma 3.10 and Theorem 3.8 of [4] are assumed to have the countable
concatenation property, then Theorem 3.10 is applicable to epi$f$ in
the proof of Theorem 3.8 of [4] and Lemma 3.17 is applicable to $V$
in the proof of Lemma 3.10 since in the two cases they automatically
have the countable concatenation property, in particular locally
$L^{0}-$convex modules currently useful in financial applications
all have such a property.

For a random locally convex module $(E,{\cal P}), \ x\in E$ and
$G\subset E$, let $d^{\ast}_{{\cal Q}}(x,G)=\wedge\{\|x-y\|_{{\cal
Q}}\,|\,y\in G\}$ for any finite subfamily ${\cal Q}$ of ${\cal P}$
and $d^{\ast}(x,G)=\vee\{d^{\ast}_{{\cal Q}}(x,G)\,|\,{\cal
Q}\subset {\cal P}$ finite $\}$. If $G$ is a nonempty ${\cal
T}_{\varepsilon,\lambda}-$closed $L^{0}-$convex subset, then we
proved in [22] that $x\in G$ iff $d^{\ast}(x,G)=0$.

{\noindent\bf Theorem 3.12([22, Theorem 3.1]).} Let $(E,{\cal P})$
be a random locally convex module over $R$ with base $(\Omega,{\cal
F},P), \ x\in E, \  G$ a nonempty ${\cal
T}_{\varepsilon,\lambda}-$closed $L^{0}-$convex subset of $E$ such
that $x \not\in G$, and $\xi$ a chosen representative of
$d^{\ast}(x,G)$. Then there exists an $f\in
E^{\ast}_{\varepsilon,\lambda}$ such that

$f(x)>\vee\{f(y)\,|\,y\in G\}$ on $[\xi>0]$, where
$[\xi>0]=\{\omega\in\Omega\,|\,\xi(\omega)>0\},\hfill{}(3.14)$

\noindent and

 $f(x)>\vee\{f(y)\,|\,y\in G\}.\hfill{}(3.14)^{\prime}$

\noindent If $R$ is replaced by $C$, then Theorem 3.12 still holds
in the following way:

$(Ref)(x)>\vee\{(Ref)(y)\,|\,y\in G\}$ on $[\xi>0]$,$\hfill{}(3.15)$

\noindent and

$(Ref)(x)>\vee\{(Ref)(y)\,|\,y\in G\}$. $\hfill{}(3.15)^{\prime}$

{\noindent\bf Theorem 3.13.} Theorem 3.12 is equivalent to Theorem
3.1 of [22].

{\noindent\bf Proof.} Since (3.14) and $(3.14)^{\prime}$ are just
(3.15) and $(3.15)^{\prime}$, respectively, in the case of a real
space, we only need to check that $(3.15)^{\prime}$ and (3.15) are
equivalent to (1) and (2) of Theorem 3.1 of [22], respectively.

Let $\xi, \  (Ref)(x)$ and $\vee\{(Ref)(y)\,|\,y\in G\}$ be the same
ones as in Theorem 3.12, and let $\eta$ and $r$ be arbitrarily
chosen representatives of $(Ref)(x)$ and $\vee\{(Ref)(y)\,|\,y\in
G\}$, respectively. Then (1) and (2) of Theorem 3.1 of [22] are
equivalent to the following (3.16) and (3.17), respectively:

$(Ref)(x)>\vee\{(Ref)(y)\,|\,y\in G\},\hfill{}(3.16)$\\
and $P([\eta>r] \ \triangle \ [\xi>0])=0,\hfill{}(3.17)$\\
where $[\eta>r]=\{\omega\in\Omega\,|\,\eta(\omega)>r(\omega)\}$, and
$[\eta>r]\triangle[\xi>0]$ denotes the symmetric difference of
$[\eta>r]$ and $[\xi>0]$.

Clearly, (3.16) is just $(3.15)^{\prime}$, and thus we only need to
prove that (3.17) is
equivalent to $(3.15)$ of Theorem 3.12.

(3.17) implies, of course, (3.15) of Theorem 3.12. On the other
hand, (3.15) of Theorem 3.12 shows that $[\eta>r]\supset[\xi>0]$,
a.s., we will prove $[\eta>r]\subset[\xi>0]$, a.s., as follows.

Otherwise, let $D=[\eta>r]\setminus[\xi>0]$, then $P(D)>0$. Thus
$I_{D}(\omega)\cdot\xi(\omega)=0$, a.s., namely $\tilde{I}_{D}\cdot
d^{*}(x,G)=0$, so that $d^{*}(\tilde{I}_{D}x,\tilde{I}_{D}G)=0$, but
from the following Lemma 3.14, we have $\tilde{I}_{D}x\in$ the
${\cal T}_{\varepsilon,\lambda}-$ closure of $\tilde{I}_{D}G$, which
means that
$\tilde{I}_{D}(Ref)(x)=(Ref)(\tilde{I}_{D}x)=\vee\{(Ref)(\tilde{I}_{D}y)\,|\,y\in
G\}=\tilde{I}_{D}\cdot(\vee\{(Ref)(y)\,|\,y\in G\})$, namely $
I_{D}(\omega)\eta(\omega)= I_{D}(\omega)r(\omega)$, a.s., which
contradicts the fact that $\eta(\omega)>r(\omega)$, a.s. on
$D$.\quad$\Box$

In the sequel of this paper, for a subset $G$ of a random locally
convex module $(E,\mathcal {P})$, $\bar G_{\varepsilon,\lambda}$
denotes the ${\cal T}_{\varepsilon,\lambda}-$ closure of $G$, and
$\bar G_{c}$ the ${\cal T}_{c}-$ closure of $G$. For any $x\in E$,
any finite subfamily $Q$ of $\mathcal {P}$ and $\varepsilon\in
L_{++}^{0}$, let $U_{Q,\varepsilon}[x]=\{y\in
E\,|\,\|x-y\|_{Q}\leq\varepsilon\}$,
$\varepsilon_{Q}^{*}(x,G)=\wedge\{\varepsilon\in L_{++}^{0}\,|\,
U_{Q,\varepsilon}[x]\,\cap\, G\neq\emptyset\}$ and
$\varepsilon^{*}(x,G)=\vee\{\varepsilon_{Q}^{*}(x,G)\,|\,Q\subset
\mathcal {P}$ finite$\}$.

{\noindent\bf Lemma 3.14}. Let $(E,\mathcal {P})$ be a random
locally convex module over $K$ with base $(\Omega,{\cal F},P)$,
$x\in E$ and $G$ a nonempty subset of $E$. then we have the
following:

(1) $d^{*}(x,G)=\varepsilon^{*}(x,G)$;

(2) $d^{*}(x,G)=d^{\ast}(x,\bar G_{\varepsilon,\lambda})=d^{\ast}(x,\bar G_{c})$.
\\If, in addition, $G$ satisfies the following:

$$\tilde{I}_{A}y+\tilde{I}_{A^{c}}z\in G  \ \textmd{for all} \  A\in{\cal F} \textmd{ and all
 }y,z\in G,\eqno(3.18),$$\\then we have the following:

 (3) $x\in\bar G_{\varepsilon,\lambda}$ iff $d^{*}(x,G)=0$.

 {\noindent\bf Proof. }\,(1). We only need to check that
 $d^{*}_{Q}(x,G)=\varepsilon^{*}_{Q}(x,G)$ for each $Q\subset\mathcal {P}$
 finite. By definition, $d^{*}_{Q}(x,G)=\wedge\{\|x-y\|_{Q}\,|\,y\in
 G\}$ and $\varepsilon^{*}_{Q}(x,G)=\wedge\{\varepsilon\in L_{++}^{0}\,|\,U_{Q,\varepsilon}[x]\,\cap\, G
 \neq\emptyset\}$. If $\varepsilon\in L_{++}^{0}$ is such that $U_{Q,\varepsilon}[x]\,\cap\, G
 \neq\emptyset$, \  namely there exists $y\in G$ such that
 $\|x-y\|_{Q}\leq\varepsilon$, which means that $\varepsilon\geq\wedge\{\|x-y\|_{Q}\,|\,y\in
 G\}=d_{Q}^{*}(x,G)$, so that $\varepsilon^{*}_{Q}(x,G)\geq
 d^{*}_{Q}(x,G)$. In the other
 direction, for any $y\in G$ and $n\in N$, it is clear that
 $\|x-y\|_{Q}\leq\|x-y\|_{Q}+\frac{1}{n}$, and
 $\|x-y\|_{Q}+\frac{1}{n}\in L_{++}^{0}$, if, let
 $\varepsilon=\|x-y\|_{Q}+\frac{1}{n}$, then we have that $y\in U_{Q,\varepsilon}[x]\,\cap\,
 G$, of course, $U_{Q,\varepsilon}[x]\,\cap\, G\neq\emptyset$, and
 thus we have that $\varepsilon=\|x-y\|_{Q}+\frac{1}{n}\geq\wedge\{\varepsilon\in L_{++}^{0}\,|\,U_{Q,\varepsilon}[x]\,\cap\, G\neq\emptyset\},\forall  \ n\in
 N$, so that $\|x-y\|_{G}\geq\varepsilon_{Q}^{*}(x,G)$, in turn $d_{Q}^{*}(x,G)=\wedge\{\|x-y\|_{Q}\,|\,y\in
 G\}\geq\varepsilon_{Q}^{*}(x,G)$.

 (2)  is clear.

 (3). Lemma 2.2 of [22] shows that $x\in F$ iff $d^{*}(x,F)=0$ for
 every ${\cal T}_{\varepsilon,\lambda}-$ closed $L^{0}-$convex
 subset $F$ of $E$, in which proof the only property (3.18) of an $L^{0}-$convex subset was used,
 and thus we have that $x\in F$ iff $d^{*}(x,F)=0$ for every ${\cal T}_{\varepsilon,\lambda}-$closed subset $F$
 with the property (3.18). It is easy to see that $\bar
 G_{\varepsilon,\lambda}$ also has the property (3.18) if $G$ does.
 Applying the result to $\bar G_{\varepsilon,\lambda}$ yields that $x\in\bar
 G_{\varepsilon,\lambda}$ iff $d^{*}(x,\bar
 G_{\varepsilon,\lambda})=0$, so that (2) has implied that $x\in\bar
 G_{\varepsilon,\lambda}$, iff $d^{*}(x,G)=0$. $\square$

 Perhaps, one would have that $x\in G$ iff $d^{*}(x,G)=0$ for every
 ${\cal T}_c-$closed subset $G$
 with the property (3.18), but it is not true since
 ${\cal T}_c$ is too strong! (see the following Lemma 3.17)
 Lemma 3.14 leads to the following equivalent variant Theorem 3.15
 of Theorem 3.12, one only needs to notice that $x\notin  \bar
 G_{\varepsilon,\lambda}$ iff $d^{*}(x,G)>0$ for every subset $G$ with
 the property (3.18), and then applies Theorem 3.12 to $x$ and $\bar G_{\varepsilon,\lambda}$.

 {\noindent\bf Theorem 3.15.} Let $(E,\cal {P})$ be a random locally convex
 module over $R$ with base $(\Omega,{\cal F},P), \  x\in
 E$ and $G$ an $L^{0}-$convex subset of $E$ such that
 $d^{*}(x,G)>0$. Then there exists an $f\in
 E_{\varepsilon,\lambda}^{*}$ such that $$f(x)>\vee\{f(y)\,|\,y\in G\} ~ \mbox{on} \ [\xi>0], \eqno(3.19)$$
 \noindent and

$$f(x)>\vee\{f(y)\,|\,y\in G\}.   \eqno(3.19)^{\prime}$$

 \noindent Where $\xi$ is an arbitrarily chosen representative of
 $d^{*}(x,G)$.

 Furthermore, if $R$ is replaced by $C$, then
 Theorem 3.15 still holds in the following way:$$(Ref)(x)>\vee\{(Ref)(y)\,|\,y\in G\} ~
 \mbox{on} ~
 [\xi>0],\eqno(3.19)^{\prime\prime}$$
 \noindent and

$$(Ref)(x)>\vee\{(Ref)(y)\,|\,y\in G\}. ~\eqno(3.19)^{\prime\prime\prime}$$

 {\noindent\bf Theorem 3.16.} Theorem 3.15 (equivalently, Theorem
 3.12) implies Theorem 3.10.

{\noindent\bf Proof. }If $x$ and $G$ satisfy the conditions of
Theorem 3.10, then
 the following Lemma 3.17 shows that $d^{*}(x,G)>0$ on $\Omega$,
 namely $[\xi>0]$ in Theorem 3.15 is just $\Omega$, so there exists
 an $f\in E_{\varepsilon,\lambda}^{*}$ such that
 $f(x)>\vee\{f(y)\,|\,y\in G\}$ on $\Omega$ by (3.19) of Theorem 3.15. Let
 $\varepsilon=\frac{1}{2}(f(x)-\vee\{f(y)\,|\,y\in G\})$, then $
 \varepsilon\in L_{++}^{0}$ satisfies the
 following:$$f(x)>\vee\{f(y)\,|\,y\in G\}+\varepsilon, \eqno(3.20)$$ $f$,
 of course, satisfies the requirement of Theorem 3.10 if one
 notices
 that $f$ is also in $E_{c}^{*}$ since
 $E_{c}^{*}=E_{\varepsilon,\lambda}^{*}$ by observing that $\cal {P}$
 has the countable concatenation property.\quad $\Box$

The following Lemma 3.17 occurred in [4, p.4015] where an outline of
its idea of proof was also given, we give its proof in detail to
find out the following Theorem 3.18.

 {\noindent\bf Lemma 3.17.} Let $(E,\cal {P})$ be a random locally
 convex module over $K$ with base
 $(\Omega,{\cal F},P), x\in E$ and $G\subset E$ a
 ${\cal T}_{c}-$closed nonempty subset such that $\hat
 I_{A}\{x\}\,\cap\,\tilde{I}_{A}G=\emptyset$ for all $A\in {\cal F}$ with
 $P(A)>0$ and such that $G$ has the countable concatenation property.
 Then $d^{*}(x,G)=\varepsilon^{*}(x,G)>0$ on $\Omega$, namely
 $\varepsilon^{*}(x,G)\wedge 1\in L_{++}^{0}$.

 {\noindent\bf Proof.} If it is not true that $d^{*}(x,G)>0$ on $\Omega$, then there
 is an $A\in {\cal F}$ such that $P(A)>0$ and $\tilde{I}_{A}\cdot
 d^{*}(x,G)=0$, so that $\tilde{I}_{A}\cdot d^{*}_{Q}(x,G)=0$ for every
 finite subfamily $Q$ of $\cal {P}$.

 We can, without loss of generality, assume that $\theta\in
 G$(otherwise, by a translation). Since $G$ has the countable
 concatenation property, it must have the property (3.18) so that
 $\{\|\tilde{I}_{A}x-\tilde{I}_{A}y\|_{Q}\,|\,y\in G\}$ is directed downwards
 for each $Q\subset \cal{P}$ finite, and it is also easy to see from
 the property (3.18) that $\tilde{I}_{A}G\subset G$, so that we can
 easily prove that $\tilde{I}_{A}G$ is also ${\cal T}_{c}-$closed.

 Now, for each fixed $Q\subset \cal{P}$ finite and each fixed
 $\alpha\in L_{++}^{0}$, we will prove that there is $y_{Q,\alpha}\in
 G$ such that $\|\tilde{I}_{A}x-\tilde{I}_{A}y_{Q,\alpha}\|_{Q}\leq\alpha$
 as follows.

 Since $\tilde{I}_{A}d^{*}_{Q}(x,G)=\wedge\{\|\tilde{I}_{A}x-\tilde{I}_{A}y\|_{Q}\,|\,y\in
 G\}=0$, then there exists a sequence $\{y_{n}\,|\,n\in N\}$ in $G$
 such that $\{\|\tilde{I}_{A}x-\tilde{I}_{A}y_n\|_{Q}\,|\,n\in N\}$ converges to $0$
 in a nonincreasing way. Let
 $\varepsilon_{n}=\|\tilde{I}_{A}x-\tilde{I}_{A}y_n\|_{Q}$ for each $n\in N$ and choose a representative
 $\varepsilon_{n}^{0}$ of $\varepsilon_{n}$ for each $n\in N$ such
 that $\varepsilon_{n}^{0}(\omega)\geq\varepsilon_{n+1}^{0}(\omega)$ for each
 $n\in N$ and $\omega\in \Omega$, and a representative $\alpha^{0}$ of
 $\alpha$ such that $[\varepsilon^{0}_{n}\leq\alpha^{0}]\uparrow\Omega$, where
 $[\varepsilon^{0}_{n}\leq\alpha^{0}]$, denoted by $E_{n}$, $=\{\omega\in
 \Omega\,|\,\varepsilon^{0}_{n}(\omega)\leq\alpha^{0}(\omega)\}$. Again let
 $A_{n}=E_{n}/E_{n-1},n\geq1$, where $E_{0}=\emptyset$, then
 $\{A_{n}\,|\,n\in N\}$ forms a countable partition of $\Omega$ to
 ${\cal F}$, so that $y_{Q,\alpha}:=\sum_{n\geq1}\tilde{I}_{A_{n}}y_{n}\in G$ and satisfies that $\|\tilde{I}_{A}x-\tilde{I}_{A}y_{Q,\alpha}\|_{Q}\leq\alpha$.

 Finally, let $\Gamma={\cal F}({\cal P})\times L_{++}^{0}=\{(Q,
 \alpha)\,|\, Q \subset {\cal P}$ finite and $\alpha\in L_{++}^{0} \}$, where ${\cal F}({\cal P})$ is the set of finite subfamilies of ${\cal P}$.
 It is easy to see that $\Gamma$ is directed upwards by the
 ordering:$(Q_{1},\alpha_{1})\leq(Q_{2},\alpha_{2})$ iff
 $Q_{1}\subset Q_{2}$ and $\alpha_{2}\leq\alpha_{1}$, so that $\{\hat
 I_{A}y_{Q,\alpha}\,,\,(Q,\alpha)\in \Gamma\}$ is a net in $\tilde{I}_{A}G$
 which is convergent to $\tilde{I}_{A}x$, and hence $\tilde{I}_{A}x\in \tilde{I}_{A}G$, which contradicts the fact that $\tilde{I}_{F}x\notin \tilde{I}_{F}G$ for all $F\in {\cal F}$ with $P(F)>0.\quad\Box$

 From the process of proof of Lemma 3.17, we can easily see that if
 $G$ has the countable concatenation property and $d^{*}(x,G)=0$ (at
 which time, take $A=\Omega$) then $x\in \bar G_{c}$, this yields a
 useful fact, namely Theorem 3.18 below, from which we have that a
 subset having the countable concatenation property has the  same
 closure under ${\cal T}_{c}$ and
 ${\cal T}_{\varepsilon,\lambda}$, in particular it is
 ${\cal T}_{c}-$closed iff it is
 ${\cal T}_{\varepsilon,\lambda}-$closed, which also derives a
 surprising fact on ${\cal T}_{c}-$ completeness, see Subsection
 3.4.

 {\noindent\bf Theorem 3.18.} Let $(E,\cal {P})$ be a random locally
 convex module, $x\in E$ and $G\subset E$ a subset having the
 countable concatenation property. Then the following are equivalent:

 (1) $x\in \bar G_{c}$;

 (2) $x\in \bar G_{\varepsilon,\lambda}$;

 (3) $d^{*}(x,G)=0$

 Theorem 3.12 (namely, Theorem 3.15) has been known for many years,
 which was mentioned in [21] without proof, where a special case of
 which was proved. We may say that Theorem 3.12 is general enough to
 meet all our needs under ${\cal T}_{\varepsilon,\lambda}$, and it
 implies Theorem 3.10 but is independent of Theorem 3.8! To
 generalize Theorems 3.8 and 3.10 to meet our further needs of
 Subsection 3.3, we present the notion of a countable concatenation
 closure (see Definition 3.19 below) to give an interesting purely
 algebraic result as follows, which provides a geometric intuition on the
 conditions
 imposed on Theorems 3.8 and 3.10, namely $\tilde{I}_{A}G\,\cap\,\tilde{I}_{A}M=\emptyset$ and $\tilde{I}_{A}\{x\}\,\cap\,\tilde{I}_{A}G=\emptyset$ for
 all $A\in {\cal F}$ with $P(A)>0$, respectively.

 {\noindent\bf Definition 3.19.} Let $E$ be a left module over the
 algebra $L^{0}({\cal F},K)$. Two countable concatenations
 $\sum_{n\geq1}\tilde{I}_{A_{n}}x_{n}$ and $\sum_{n\geq1}\tilde{I}_{B_{n}}y_{n}$ are called equal if $\tilde{I}_{A_{i}\,\cap\,
 B_{j}}x_{i}=\tilde{I}_{A_{i}\,\cap\, B_{j}}y_{j},\forall i,j\in N$. For any
 subset $G$ of $E$, the set $\mathcal {C}_{cc}(G)=\{\sum_{n\geq1}\tilde{I}_{A_{n}}x_{n}|\sum_{n\geq1}\tilde{I}_{A_{n}}x_{n}$ is a countable
 concatenation and each $x_{n}\in G \}$ is called the countable
 concatenation closure of $G$.

 {\noindent\bf Theorem 3.20.} Let $E$ be a left module over the
 algebra $L^{0}({\cal F},K)$, $M$ and $G$ any two nonempty subsets
 of $E$ such that $\tilde{I}_{A}M+\tilde{I}_{A^{c}}M\subset M$ and $\tilde{I}_{A}G+\tilde{I}_{A^{c}}G\subset G$. If $\mathcal
 {C}_{cc}(M)\,\cap\,\mathcal {C}_{cc}(G)=\emptyset$, then there exists an
 ${\cal F}-$measurable subset $H(M,G)$ unique a.s. such that the
 following are satisfied:

 (1) $P(H(M,G))>0$;

 (2) $\tilde{I}_{A}M\,\cap\,\tilde{I}_{A}G=\emptyset$ for all $A\in
 {\cal F},A\subset H(M,G)$ with $P(A)>0$;

 (3) $\tilde{I}_{A}M\,\cap\,\tilde{I}_{A}G\neq\emptyset$ for all $A\in
 {\cal F},A\subset \Omega\backslash H(M,G)$ with $P(A)>0$.

{\noindent\bf Proof.} Let ${\cal E}=\{A\in {\cal
F}\,|\,\tilde{I}_{A}M\,\cap\,\tilde{I}_{A}G\neq\emptyset\}$. Then
${\cal E}$ is directed upwards: in
 fact, for any $A$ and $B\in{\cal E}$ there exist $x_{1},x_{2}\in
 M$ and $y_{1},y_{2}\in G$ such that $\tilde{I}_{A}x_{1}=\hat
 I_{A}y_{1}$ and $\tilde{I}_{B}x_{2}=\tilde{I}_{B}y_{2}$. Since $M$ and $G$
 are nonempty, take $x_{0}$ in $M$ and $y_{0}$ in $G$, and let
 $x=\tilde{I}_{A}x_{1}+\tilde{I}_{B\backslash A}x_{2}+\tilde{I}_{(A\cup
 B)^{c}}x_{0}$ and $y=\tilde{I}_{A}y_{1}+\tilde{I}_{B\backslash
 A}y_{2}+\tilde{I}_{(A\cup B)^{c}}y_{0}$, then $\tilde{I}_{A\cup B}x=\tilde{I}_{A\cup B}y\in \tilde{I}_{A\cup B}M\,\cap\,\tilde{I}_{A\cup B}G$ by noticing
 $x\in M$ and $y\in G$.

 Define $H(M,G)=\Omega\backslash esssup({\cal E})$, then $H(M,G)$,
 obviously,satisfies (2) and (3). We will verify that $H(M,G)$ also
 has the property (1). In fact, if $P(H(M,G))=0$, then
 $P(esssup({\cal E}))=1$, let $\{D_{n}\,|\,n\in N\}$ be a
 nondecreasing sequence of ${\cal E}$ such that
 $D_{n}\uparrow\Omega$, then there exist two sequences $\{x_{n}\,|\,n\in
 N\}$ in $M$ and $\{y_{n}\,|\,n\in N\}$ in $G$ such that $\tilde{I}_{D_{n}}x_{n}=\tilde{I}_{D_{n}}y_{n},\forall n\in N$. Let
 $A_{n}=D_{n}\setminus D_{n-1},\forall n\geq1$, where
 $D_{0}=\emptyset$, then $\sum_{n\geq1}\tilde{I}_{A_{n}}x_{n}=\sum_{n\geq1}\tilde{I}_{A_{n}}y_{n}=\mathcal
 {C}_{cc}(M)\,\cap\,\mathcal {C}_{cc}(G)$, which is a
 contradiction.\quad $\Box$

 {\noindent\bf Definition 3.21.} Let $E$,$M$ and $G$ be the same as in
 Theorem 3.20 such that $\mathcal {C}_{cc}(M)\,\cap\,\mathcal
 {C}_{cc}(G)=\emptyset$, then $H(M,G)$ is called the hereditarily
 disjoint stratification of $H$ and $M$, and $P(H(M,G))$ is called the
 hereditarily disjoint probability of $H$ and $G$.

 \subsection*{3.2. The hereditarily disjoint probability and more general forms of hyperplane separation theorems}
 \quad\, First, we state the main result of this
 subsection---Theorems 3.22 and 3.23, whose proofs follow from Lemma
 3.24.

 {\noindent\bf Theorem 3.22.} Let $(E,\cal {P})$ be a random locally
 convex module over $R$ with base
 $(\Omega,{\cal F},P)$, $M$ and $G$ two nonempty $L^{0}-$convex
 subsets such that $G$ is ${\cal T}_c-$open and $\mathcal {C}_{cc}(G)\,\cap\,\mathcal
 {C}_{cc}(M)=\emptyset$. Then there exists an $f\in E_{c}^{*}$ such
 that $$f(x)<f(y) \textmd{ on } H(M,G) \textmd{ for all } x\in G \mbox{ and } y\in M, \eqno(3.21)$$
 \noindent and
 $$f(x)<f(y)  \textmd{ for all } x\in G \mbox{ and } y\in M. \eqno(3.21)^{\prime}$$

 If $R$ is replaced by $C$, then Theorem 3.22
 still holds in the following way: $$(Ref)(x)<(Ref)(y) \textmd{ on } H(M,G) \textmd{ for all } x\in G \mbox{ and } y\in M, \eqno(3.22)$$
 \noindent and
 $$(Ref)(x)<(Ref)(y) \textmd{ for all } x\in G \mbox{ and } y\in M. \eqno(3.22)^{\prime}$$

In the following Theorem 3.23, we only need to notice that
$C_{cc}(G)=G$ and $H(\{x\},G)$
 is just $[\xi>0]$ in Theorem 3.15.

 {\noindent\bf Theorem 3.23.} Let $(E,\cal {P})$ be a random locally
 convex module over $R$ with base
 $(\Omega,{\cal F},P)$, such that $\mathcal {P}$ has the countable
 concatenation property, $x\in E$ and $G$ a nonempty
 ${\cal T}_{c}-$closed $L^{0}-$convex subsets of $E$ such that
 $x\notin G$ and $G$ has the countable concatenation property. Then
 there exists an $f\in E_{c}^{*}$ such that  $$f(x)>\vee\{f(y)\,|\,y\in G\} \mbox{
 on }
 H(\{x\},G), \eqno(3.23)$$
 \noindent and
$$f(x)>\vee\{f(y)\,|\,y\in G\}. \eqno(3.23)^{\prime}$$

  If $R$ is replaced by
 $C$, then Theorem 3.23 still holds in the following
 way:$$(Ref)(x)>\vee\{(Ref)(y)\,|\,y\in G\} \mbox{
 on } H(\{x\},G), \eqno(3.24)$$
 \noindent and
 $$(Ref)(x)>\vee\{(Ref)(y)\,|\,y\in G\}. \eqno(3.24)^{\prime}$$

 {\noindent\bf Lemma 3.24.} Let $(E,\cal {P})$ be a random locally
 convex module over $K$ with base
 $(\Omega,{\cal F},P)$, $M$ a ${\cal T}_{c}-$closed subset of
 $E$ such that $\tilde{I}_{A}M+\tilde{I}_{A^{c}}M\subset M$, for all $A\in
 {\cal F}$, and $G$ a ${\cal T}_{c}-$open subset of $E$ such
 that $\tilde{I}_{A}G+\tilde{I}_{A^{c}}G\subset G$, for all $A\in
 {\cal F}$. Then for each $A\in {\cal F}$ with $P(A)>0$, $\tilde{I}_{A}M$
 is relatively ${\cal T}_{c}-$closed in $\tilde{I}_{A}E$ and $\tilde{I}_{A}G$ is relatively ${\cal T}_{c}-$open in $\tilde{I}_{A}E$.

 {\noindent\bf Proof.} We can assume that $\theta\in G$ and $\theta\in M$ (otherwise,by a
 translation), respectively, then $\tilde{I}_{A}G\subset G$, and $\tilde{I}_{A}M\subset
 M$, so that $\tilde{I}_{A}E\,\cap\, M=\tilde{I}_{A}M$ and $\tilde{I}_{A}E\,\cap\, G=\tilde{I}_{A}G.\quad\Box$

We can now prove Theorem 3.22.

 {\noindent\bf Proof of Theorem 3.22.} Let $\Omega'=H(M,G)$, ${\cal F}'=\Omega'\,\cap\,{\cal F}=\{\Omega'\,\cap\, F\,|\,F\in
 {\cal F}\}$ and $P':{\cal F}'\rightarrow [0,1]$ be defined
 by $P'(\Omega'\,\cap\, F)=P(\Omega'\,\cap\, F)/P(\Omega'),\forall F\in
 {\cal F}$. Take $E'=\tilde{I}_{\Omega'}E$, ${\cal P}'=\{\|\cdot\|_{E'}\,|\,\|\cdot\|\in
 \cal{P}\}$, $M'=\tilde{I}_{\Omega'}M$, $G'=\tilde{I}_{\Omega'}G$ and
 consider $(E',\cal {P}')$ as a random locally convex module with
 base $(\Omega',{\cal F}',P')$. Then $M'$ and $G'$ satisfy the
 condition of Theorem 3.8, so that there exists an $f'\in (E')_{c}^{*}$ such that $$f'(x)<f'(y)\mbox{ on } \Omega'\mbox{ for all } x\in G'\mbox{ and } y\in M'. \eqno(3.25)$$

 By Theorem 2.9 $f'$ has an extension $f''$
$\in E_{c}^{\ast}$. Now let $f=\tilde{I}_{H(M,G)}f''$, then
$f(x)=0,\forall x\in \tilde{I}_{H(M,G)^{c}}E$ and
$f(x)=f'(x),\forall x\in \tilde{I}_{H(M,G)}E$, so that $f$ satisfies
all the requirements of Theorem 3.22. \quad$\Box$

{\noindent\bf Proof of Theorem 3.23.} It is completely similar to
the proof of Theorem 3.22 (In fact, it can also be derived directly
from Theorem 3.15).\quad$\Box$

\subsection*{3.3. Closed $L^{0}-$convex subsets with the countable concatenation property}

{\noindent\bf Definition 3.25.} Let $(E,{\cal P})$ be a random
locally convex module over $K$ with base $(\Omega,{\cal F},P)$. For
each $f\in E_{c}^{\ast}, |f(\cdot)|:E\rightarrow L^{0}_{+}$ is
clearly an $L^{0}-$seminorm, so that $(E,\{|f(\cdot)||f\in
E_{c}^{\ast}\})$ is a random locally convex module, whose locally
$L^{0}-$convex topology, denoted by $\sigma_c(E,E_{c}^{\ast})$,
called the weak locally $L^{0}-$convex topology of $E$. Similarly,
we may have the weak $(\varepsilon,\lambda)-$topology
$\sigma_{\varepsilon,\lambda}(E,E_{\varepsilon,\lambda}^{\ast})$ of
$E$. In particular when
$E_{\varepsilon,\lambda}^{\ast}=E_{c}^{\ast}$, we briefly write
$\sigma_c(E,E^{\ast})$ and
$\sigma_{\varepsilon,\lambda}(E,E^{\ast})$ for
$\sigma_c(E,E_{c}^{\ast})$ and
$\sigma_{\varepsilon,\lambda}(E,E_{\varepsilon,\lambda}^{\ast})$,
respectively.

{\noindent\bf Remark 3.26.} For a random locally convex module
$(E,\cal P)$, dually, we may have the weak-star locally
$L^{0}-$convex topology $\sigma_c(E_{c}^{\ast},E)$ of
$E_{c}^{\ast}$, and the weak-star $(\varepsilon,\lambda)-$topology
$\sigma_{\varepsilon,\lambda}(E_{\varepsilon,\lambda}^{\ast},E)$ of
$E_{\varepsilon,\lambda}^{\ast}$, which can be briefly denoted by
$\sigma_c(E^{\ast},E)$ and
$\sigma_{\varepsilon,\lambda}(E^{\ast},E)$ when
$E_{c}^{\ast}=E_{\varepsilon,\lambda}^{\ast}$, respectively.

The main result of the subsection is the following:

{\noindent\bf Theorem 3.27.} Let $(E,\cal P)$ be a random locally
convex module over $K$ with base $(\Omega,{\cal F},P)$ and $G$ an
$L^{0}-$convex subset of $E$ such that $G$ has the countable
concatenation property. Then we have the following equivalent
statements:

(1) $G$ is ${\cal T}_{c}-$closed;

(2) $G$ is ${\cal T}_{\varepsilon,\lambda}-$closed;

(3) $G$ is
$\sigma_{\varepsilon,\lambda}(E,E_{\varepsilon,\lambda}^{\ast})-$closed.\\
Further, if $\cal P$ has the countable concatenation property, then
the above three are equivalent to the following:

(4) $G$ is $\sigma_c(E,E_{c}^{\ast})-$closed.

{\noindent\bf Proof.} $(1)\Leftrightarrow(2)$ has been proved in
Theorem 3.18 and $(2)\Leftrightarrow(3)$ has been proved in
Corollary 3.4 of [22] for any $L^{0}-$convex subset $G$.

If $\cal P$ has the countable concatenation property, completely
similar to the proof of the classical Mazur's theorem, cf. [3], one
can see $(1)\Leftrightarrow(4)$ by $(3.24)^{\prime}$ of Theorem
3.23. $\Box$

\subsection*{3.4 Completeness}
\quad\, The main result of this subsection is the following:

{\noindent\bf Theorem 3.28.} Let $(E,\cal P)$ be a random locally
convex module over $K$ with base $(\Omega,{\cal F},P)$. Then $E$ is
${\cal T}_{c}-$complete if $E$ is ${\cal
T}_{\varepsilon,\lambda}-$complete. Furthermore, if $E$ is ${\cal
T}_{c}-$complete and has the countable concatenation property, then
$E$ is also ${\cal T}_{\varepsilon,\lambda}-$complete.

The first part of Theorem 3.18 is easy as pointed out in Section
1.4. However, the proof of the second part of it is a delicate
matter since ${\cal T}_c$ is much stronger than ${\cal
T}_{\varepsilon,\lambda}$, and it seems to be not an easy work for
us to give a direct proof even for the case of an $RN$ module.
However, we may give a clever proof of it by using Theorem 3.18 and
the following Lemma 3.29. First, let us recall some terminologies as
follows.

Let  $(E,\cal P)$ be a random locally convex module over $K$ with
base $(\Omega,{\cal F},P)$. Two ${\cal
T}_{\varepsilon,\lambda}-$Cauchy nets
$\{x_{\alpha},\alpha\in\Gamma\}$ and $\{y_{\beta},\beta\in\Lambda\}$
in $E$ are called equivalent if
$\{\|x_{\alpha}-y_{\beta}\|,(\alpha,\beta)\in\Gamma\times\Lambda\}$
converges to 0 in probability $P$ for each $\|\cdot\|\in\cal P$.
Since the sum of $\{x_{\alpha},\alpha\in\Gamma\}$ and
$\{y_{\beta},\beta\in\Lambda\}$ is defined, as usual, to be
$\{x_{\alpha}+y_{\beta},(\alpha,\beta)\in\Gamma\times\Lambda\}$,
which motivates us to do the following thing: If $E$ has the
countable concatenation property, $\{A_n\,|\,n\in N\}$ is a
countable partition of $\Omega$ to ${\cal F}$ and
$\{\{x_{\alpha_n},\alpha_n\in\Gamma_n\}\,|\,n\in N\}$ is a sequence
of Cauchy nets in $E$, then we naturally define their countable
concatenation $\sum_{n\geq
1}\tilde{I}_{A_n}\{x_{\alpha_n},\alpha_n\in\Gamma_n\}=\{\sum_{n\geq
1}\tilde{I}_{A_n}x_{\alpha_n},(\alpha_1,\alpha_2,\cdots,\alpha_n,\cdots)\in\prod_{n\geq
1}\Gamma_n\}$, and it is easy to check that this is again a ${\cal
T}_{\varepsilon,\lambda}-$Cauchy net by noticing that
$P(A_n)\rightarrow 0$ as $n\rightarrow\infty$, so that this
countable concatenation is well defined. According to the same fact
that $P(A_n)\rightarrow 0$ as $n\rightarrow\infty$, we can verify
that if $\{y_{\beta_n}:\beta_n\in\Lambda_n\}$ is equivalent to
$\{x_{\alpha_n},\alpha_n\in\Gamma_n\}$ for each $n\in N$, then
$\{\sum_{n\geq
1}\tilde{I}_{A_n}x_{\alpha_n},(\alpha_1,\alpha_2,\cdots,\alpha_n,\cdots)\in\prod_{n\geq
1}\Gamma_n\}$ is still equivalent to $\{\sum_{n\geq
1}\tilde{I}_{A_n}y_{\beta_n},(\beta_1,\beta_2,\cdots,\beta_n,\cdots)\in\prod_{n\geq
1}\Lambda_n\}$. These observations lead to Lemma 3.29 below.

{\noindent\bf Lemma 3.29.} Let $(E,\cal P)$ be a random locally
convex module over $K$ with base $(\Omega,{\cal F},P)$ such that $E$
has the countable concatenation property. For a ${\cal
T}_{\varepsilon,\lambda}-$Cauchy net
$\{x_{\alpha},\alpha\in\Gamma\},[\{x_{\alpha},\alpha\in\Gamma\}]$
denotes its ${\cal T}_{\varepsilon,\lambda}-$equivalence class, the
${\cal T}_{\varepsilon,\lambda}-$equivalence class of a constant net
with value $x\in E$ is denoted by $[x]$. Let
$\tilde{E}_{\varepsilon,\lambda}=\{[\{x_{\alpha},\alpha\in\Gamma\}]\,|\,\{x_{\alpha},\alpha\in\Gamma\}$
is a ${\cal T}_{\varepsilon,\lambda}-$Cauchy net$\}$. The module
operations are defined as follows:

$$[\{x_{\alpha},\alpha\in\Gamma\}]+[\{y_{\beta},\beta\in\Lambda\}]:=[\{x_{\alpha}+y_{\beta},(\alpha,\beta)\in\Gamma\times\Lambda\}],$$
$$\xi[\{x_{\alpha},\alpha\in\Gamma\}]=[\{\xi
x_{\alpha},\alpha\in\Gamma\}].$$

Further, each $\|\cdot\|\in\cal P$ induces an $L^{0}-$seminorm on
$\tilde{E}_{\varepsilon,\lambda}$, still denoted by $\|\cdot\|$, so
that $\|[\{x_{\alpha},\alpha\in\Gamma\}]\|$=the limit of convergence
in probability $P$ of $\{\|x_{\alpha}\|,\alpha\in\Gamma\}$.

Then $(\tilde{E}_{\varepsilon,\lambda},\cal P)$ is a ${\cal
T}_{\varepsilon,\lambda}-$complete random locally convex module over
$K$ with base $(\Omega,{\cal F},P)$ such that
$\tilde{E}_{\varepsilon,\lambda}$ still has the countable
concatenation property, called the ${\cal
T}_{\varepsilon,\lambda}-$completion of $(E,{\cal P})$, further
$(E,{\cal P})$ is ${\cal P}-$isometrically isomorphic with a dense
submodule $\{[x]|x\in E\}$ of $\tilde{E}_{\varepsilon,\lambda}$.

{\noindent\bf Proof.} From the proof of completion of a linear
topological space, we can first have that
$(\tilde{E}_{\varepsilon,\lambda},{\cal P})$ is ${\cal
T}_{\varepsilon,\lambda}-$complete. Further, let
$\{\{x_{\alpha_n},\alpha_n\in\Gamma_n\}\,|\,n\in N\}$ be a sequence
of ${\cal T}_{\varepsilon,\lambda}-$Cauchy nets in $E$ and define
$\pi_{\alpha}^{n}=x_{\alpha_n},\forall\alpha\in\prod_{n\geq
1}\Gamma_n$ and $n\in N$, then
$\{x_{\alpha_n},\alpha_n\in\Gamma_n\}$ is equivalent to
$\{\pi_{\alpha}^{n}, \alpha\in\prod_{n\geq 1}\Gamma_n\}$ for each
$n\in N$.

Thus for a countable partition $\{A_n\,|\,n\in N\}$ of $\Omega$ to
${\cal F}$ we have that $\tilde{I}_{A_m}[\{\sum_{n\geq
1}\tilde{I}_{A_n}x_{\alpha_n}, (\alpha_1,$ $\alpha_2, \cdots,
\alpha_n, \cdots)\in\prod_{n\geq
1}\Gamma_n\}]=[\{\tilde{I}_{A_m}\pi_{\alpha}^{m},\alpha\in\prod_{n\geq
1}\Gamma_n\}]$(by the definition of the module
multiplication)=$\tilde{I}_{A_m}[\{\pi_{\alpha}^{m},\alpha\in\prod_{n\geq
1}\Gamma_n\}]=\tilde{I}_{A_m}[\{x_{\alpha_m},\alpha_m\in\Gamma_m\}]$
for each $m\in N$, so that $\sum_{m\geq
1}\tilde{I}_{A_m}[\{x_{\alpha_m},\alpha_m\in\Gamma_m\}]=[\{\sum_{n\geq
1}\tilde{I}_{A_n}x_{\alpha_n},
(\alpha_1,\alpha_2,\cdots,\alpha_n,\cdots)\in\prod_{n\geq
1}\Gamma_n\}]\in \tilde{E}_{\varepsilon,\lambda}$, namely
$\tilde{E}_{\varepsilon,\lambda}$ has the countable concatenation
property. \quad $\Box$

{\noindent\bf Remark 3.30.} Although every random locally convex
module $(E,{\cal P})$ admits a ${\cal T}_c-$completion
$(\tilde{E}_c,{\cal P})$ and a ${\cal
T}_{\varepsilon,\lambda}-$completion
$(\tilde{E}_{\varepsilon,\lambda},{\cal P})$ such that $(E,{\cal
P})$ is $\cal P-$isometrically isomorphic with a dense submodule of
either of the latter two, but since ${\cal T}_c$ is so strong that a
countable concatenation of a sequence of ${\cal T}_c-$Cauchy nets is
not necessarily well defined, so that we can not give Lemma 3.29 for
${\cal T}_c-$topology. Lemma 3.29 is necessary since the countable
concatenation property of $E$ is reserved in
$\tilde{E}_{\varepsilon,\lambda}$ in a proper way when $E$ and
$\{[x]\,|\,x\in E\}$ are identified.

We can now prove Theorem 3.28.

{\noindent\bf Proof of Theorem 3.28.} Let
$(\tilde{E}_{\varepsilon,\lambda},{\cal P})$ be the ${\cal
T}_{\varepsilon,\lambda}-$completion of $(E,{\cal P})$ as in Lemma
3.29 and regard $E$ as a subset of
$\tilde{E}_{\varepsilon,\lambda}$, then Theorem 3.18 shows that
$\bar{E}_c=\bar{E}_{\varepsilon,\lambda}=\tilde{E}_{\varepsilon,\lambda}$.
On the other hand, since $E$ is ${\cal T}_c-$complete, we always
have that $\bar{E}_c=E$, so that
$\tilde{E}_{\varepsilon,\lambda}=E$, namely $E$ must be ${\cal
T}_{\varepsilon,\lambda}-$complete.\quad $\Box$

{\noindent\bf Remark 3.31.} Theorem 3.28 is a powerful result, for
example, Kupper and Vogelpoth proved in [25] that $L^{p}_{{\cal
F}}({\cal E})$ and $L^{\varphi}_{{\cal F}}({\cal E})$ are ${\cal
T}_c-$complete, then they must be ${\cal
T}_{\varepsilon,\lambda}-$complete by the second part of Theorem
3.28, since they both have the countable concatenation property. On
the other hand, it is easy to verify that $L^{1}_{{\cal F}}({\cal
E})$ is ${\cal T}_c-$complete, as to $L^{p}_{{\cal F}}({\cal E})$
when $p>1$ they are, obviously, ${\cal T}_c-$complete by the first
part of Theorem 3.28 since they are all the random conjugate spaces
of some $RN$ modules and ${\cal T}_{\varepsilon,\lambda}-$complete.

\section*{4. The theory of random conjugate spaces of random normed
modules under the locally $L^{0}-$convex topology}

\quad\, Since the $(\varepsilon,\lambda)-$topology is rarely a locally
convex topology in the sense of traditional functional analysis,
consequently, the theory of traditional conjugate spaces universally
fails to serve for the deep development of $RN$ modules under the
$(\varepsilon,\lambda)-$topology, see [22] for details. It is under
such a background that the theory of random conjugate spaces of
random normed modules has been developed and has been being centered
at our previous work, and in fact it is also the most difficult and
deepest part of our previous work, cf. [8, 10, 16, 19, 23].

The locally $L^{0}-$convex topology has the nice convexity, but it
is too strong, it is, certainly, also rather difficult to establish
the corresponding results of [8, 10, 16, 19, 23] under the locally
$L^{0}-$convex topology in a direct way. Considering that an $RN$
module has the same random conjugate space under the two kinds of
topologies, even many results are independent of a special choice of
the two kinds of topologies, we can now establish the corresponding
${\cal T}_{c}-$variants of those deep results previously established
in [8, 10, 16, 19, 23] under ${\cal T}_{\varepsilon,\lambda}$, and
based on Section 3.4 this has become an easy matter in an indirect
manner!

It should also be pointed out that there are many results in Section
4 in which the hypothesis ``Let $(E,\|\cdot\|)$ be ${\cal
T}_{c}-$complete and have the countable concatenation property''
occurs, the hypothesis automatically reduces to ``Let
$(E,\|\cdot\|)$ be ${\cal T}_{\varepsilon,\lambda}-$complete'' in
their ${\cal T}_{\varepsilon,\lambda}-$prototypes, since ${\cal
T}_{\varepsilon,\lambda}-$completeness has implied the countable
concatenation property. Besides, the reader should bear in mind that
$E^{\ast}_{c}=E^{\ast}_{\varepsilon,\lambda}$, denoted by
$E^{\ast}$, for an $RN$ module.
\subsection*{4.1. Riesz's representation theorems and the important connection between random conjugate spaces and classical conjugate spaces}

\quad\, In this subsection, we will give the Riesz's type of
representation theorems of random conjugate spaces of three
extensive classes of $RN$ modules. The main results are Theorem 4.3,
Theorem 4.4 and Theorem 4.8.

To give the first of Riesz's type of representation theorems, let us
first recall the notion of a random inner product module, if we do
not intend to mention the notion of a random inner product space,
then the notion of a random inner product module introduced in [11]
and already employed in [23] is exactly the following:

{\noindent\bf Definition 4.1([11]).} An ordered pair
$(E,\langle\cdot,\cdot\rangle)$ is called {\it{a random inner
product module}} (briefly, an $RIP$ module) over $K$ with base
$(\Omega,{\cal F},P)$ if $E$ is a left module over the algebra
$L^{0}({\cal F},K)$ and $\langle\cdot,\cdot\rangle$ is a mapping
from $E\times E$ to $L^{0}({\cal F},K)$ such that the following are
satisfied:

(1) $\langle x,x\rangle\in L^{0}_+$, and $\langle x,x\rangle=0$ iff
$x=\theta$ (the null vector of $E$);

(2) $\langle x,y\rangle=\overline{\langle y,x\rangle},\forall x,y\in
E$ where $\langle y,x\rangle$ denotes the complex conjugate of
$\langle y,x\rangle$;

(3) $\langle\xi x,y\rangle=\xi\langle x,y\rangle,\forall\xi\in
L^{0}({\cal F},K)$, and $x,y\in E$;

(4) $\langle x+y,z\rangle=\langle x,z\rangle+\langle
y,z\rangle,\forall x,y,z\in E$.\\
Where $\langle x,y\rangle$ is called the random inner product
between $x$ and $y$; If $\langle x,y\rangle=0$, then $x$ and $y$ are
called {\it{orthogonal}}, denoted by $x\perp y$, furthermore
$M^{\perp}=\{y\in E\,|\,\langle x,y\rangle=0,\forall x\in M\}$ is
called the {\it{orthogonal complement}} of $M$.

The Schwartz inequality: $|\langle
x,y\rangle|\leqslant\|x\|\cdot\|y\|,\forall x,y\in E$, was proved in
[11], where $\|\cdot\|:E\rightarrow L^{0}_+$ defined by
$\|x\|=\sqrt{\langle x,x\rangle},\forall x\in E$, is thus an
$L^{0}-$norm such that $(E,\|\cdot\|)$ is an $RN$ module, called the
$RN$ module derived form $(E,\langle\cdot,\cdot\rangle)$

Let us first recall from [13] and its references there the notions
of random elements and random variables.

{\noindent\bf Example 4.2.} Let $(H,\langle\cdot,\cdot\rangle)$ be a
Hilbert space over $K$. $A$ mapping $V:\Omega\rightarrow H$ is
called {\it{an ${\cal F}-$random element}} if
$V^{-1}(B):=\{\omega\in\Omega\,|\,V(\omega)\in B\}\in{\cal F}$ for
each open subset $B$ of $H$, further an ${\cal F}-$random element
with its range finite is called {\it{an ${\cal F}-$simple random
element}}. A mapping $V:\Omega\rightarrow H$ is called an{\it{
${\cal F}-$random variable}} if there exists a sequence
$\{V_{n}\,|\,n\in N\}$ of ${\cal F}-$simple random elements such
that $\|V_{n}(\omega)-V(\omega)\|\rightarrow 0$ as $n\rightarrow
\infty$ for each $\omega\in\Omega$. Denote by $L^{0}({\cal F},H)$
the linear space of equivalence classes of $H-$valued ${\cal
F}-$random variables on $\Omega$, which is a left module over the
algebra $L^{0}({\cal F},K)$ under the module multiplication $\xi
x:=$ the equivalence class of $\xi^{0}\cdot x^{0}$ defined by
$(\xi^{0}\cdot x^{0})(\omega)=\xi^{0}(\omega)\cdot
x^{0}(\omega),\forall\omega\in\Omega$, where $\xi^{0}$ and $x^{0}$
are arbitrarily chosen representatives of $\xi\in L^{0}({\cal F},K)$
and $x\in L^{0}({\cal F},H)$.

$L^{0}({\cal F},H)$ becomes an $RIP$ module
over $K$ with base $(\Omega,{\cal F},P)$ under the random inner
product induced from the inner product $\langle\cdot,\cdot\rangle$,
still denoted by $\langle\cdot,\cdot\rangle$. Namely, for any
$x,y\in L^{0}({\cal F},H)$ with respective representatives
$x^{0},y^{0}$, we have $\langle x,y\rangle=$ the equivalence class
of $\langle x^{0},y^{0}\rangle$ defined by $\langle
x^{0},y^{0}\rangle(\omega)=\langle
x^{0}(\omega),y^{0}(\omega)\rangle,\forall\omega\in\Omega$.

It is easy to see that $(L^{0}({\cal F},H),\|\cdot\|)$ is ${\cal
T}_{\varepsilon,\lambda}-$complete, so that $L^{0}({\cal F},H)$ is
${\cal T}_{c}-$complete and obviously has the countable
concatenation property.

{\noindent\bf Theorem 4.3.}  Let $(E,\langle\cdot,\cdot\rangle)$
be a ${\cal T}_{c}-$complete $RIP$ module over $K$ with base
$(\Omega,{\cal F},P)$ such that $E$ has the countable
concatenation property. Then for every $f\in E^{\ast}$ there
exists a unique $\pi(f)\in E$ such that $f(x)=\langle
x,\pi(f)\rangle,\forall x\in E$ and such that
$\|f\|^{\ast}=\|\pi(f)\|$. Finally the induced mapping
$\pi:E^{\ast}\rightarrow E$ is a surjective conjugate isomorphism,
namely $\pi(\xi f+\eta
g)=\bar{\xi}\pi(f)+\bar{\eta}\pi(g),\forall\xi,\eta\in L^{0}({\cal
F},K)$ and $f,g\in E^{\ast}$.

{\noindent\bf Proof.} By Theorem 3.28, we have that $E$ is ${\cal
T}_{\varepsilon,\lambda}-$complete. It follows immediately that the
main result of [23] is just what we want to prove. $\Box$

Before giving Theorem 4.4, let us first recall two important
examples of $RN$ modules: Let $L^{0}({\cal F},B)$ the $RN$ module of
equivalence classes of random variables from $(\Omega,{\cal F},P)$
to a normed space $(B,\|\cdot\|)$ over $K$, its construction is
similar to Example 4.2, also see [19] for details.

Let $B'$ be the classical conjugate space of $B$. Then a mapping
$q:\Omega\rightarrow B'$ is called a $\mbox{w}^{\ast}-$random
variable if the composite function $\langle b,q\rangle$ defined by
$\langle b,q\rangle(\omega)=\langle
b,q(\omega)\rangle,\forall\omega\in\Omega$, is a {\it{$K-$valued
random variable}} for each fixed $b\in B$, where $\langle\cdot,
\cdot\rangle:B\times B^{\prime}\rightarrow K$ denotes the natural
pairing between $B$ and $B'$. Two $\mbox{w}^{\ast}-$random
variables $q_1$ and $q_2$ are called $\mbox{w}^{\ast}-$equivalent
if $\langle b,q_1\rangle$ and $\langle b,q_2\rangle$ are
equivalent for each fixed $b\in B$. For each
$\mbox{w}^{\ast}-$random variable $q$, since $|\langle
b,q(\omega)\rangle|\leq\|q(\omega)\|$ for each $\omega\in\Omega$
and $b\in B$ such that $\|b\|\leq1$, esssup $(\{|\langle
b,q\rangle|\,|\,b\in B$ and $\|b\|\leq1\}$) is a nonnegative
real-valued random varable.

Let $L^{0}({\cal F},B',\mbox{w}^{\ast})$ be the linear space of
$\mbox{w}^{\ast}-$equivalence classes of $B'-$valued
$\mbox{w}^{\ast}-$random variables on $\Omega$. Like $L^{0}({\cal
F},H)$ in Example 4.2, $L^{0}({\cal F},B',\mbox{w}^{\ast})$ can
naturally becomes a left module over the algebra $L^{0}({\cal
F},K)$. Finally, for each $x\in L^{0}({\cal F},B',\mbox{w}^{\ast})$,
if we define its random norm $\|x\|$ by $\|x\|=$ the equivalence
class of esssup $(\{|\langle b,x^{0}\rangle||\,b\in B$ and
$\|b\|\leq1\})$, where $x^{0}$ is a representative of $x$, then
$L^{0}({\cal F},B',\mbox{w}^{\ast})$ is an $RN$ module over $K$ with
base $(\Omega,{\cal F},P)$. Finally, for any $x\in L^{0}({\cal F},
B)$ and $y\in L^{0}({\cal F}, B', \mbox{w}^{\ast})$, as usual, the
natural pairing $\langle x,y\rangle$ between $x$ and $y$ can be
defined as the equivalence class of the natural pairing between
their respective representatives.

{\noindent\bf Theorem 4.4.} Let $(\Omega,{\cal F},P)$ be a complete
probability space. Then $L^{0}({\cal F},B',\mbox{w}^{\ast})$ is
isomorphic with the random conjugate space of $L^{0}({\cal F},B)$,
denoted by $(L^{0}({\cal F},B))^{\ast}$, in a random-norm preserving
manner under the canonical mapping $T: L^{0}({\cal
F},B',\mbox{w}^{\ast})\rightarrow(L^{0}({\cal F},B))^{\ast}$ defined
as follows. For each $f\in L^{0}({\cal F},B',\mbox{w}^{\ast}),T_f$,
denoting $T(f),:L^{0}({\cal F},B)\rightarrow L^{0}({\cal F},K)$ is
given by $T_f(g)=\langle g,f\rangle,\forall g\in L^{0}({\cal F},B)$.
Further, if $L^{0}({\cal F},B',\mbox{w}^{\ast})$ is replaced by
$L^{0}({\cal F},B')$, then a sufficient and necessary condition for
$T:L^{0}({\cal F},B')\rightarrow(L^{0}({\cal F},B))^{\ast}$ to be
again an isometric isomorphism is that $B'$ has the
Randon-Nikod\'{y}m property with respect to $(\Omega,{\cal F},P)$.

{\noindent\bf{Proof.}} Since $(\Omega,{\cal F},P)$ is complete,
$L^{0}({\cal F},B',\mbox{w}^{\ast}),L^{0}({\cal F},B)$ and
$L^{0}({\cal F},B')$ are equivalent to
$L(P,B',\mbox{w}^{\ast}),L(P,B)$ and $L(P,B')$ in [8], so that our
desired results follow immediately from the main results of
[8].\quad $\Box$

{\noindent\bf Remark 4.5.} Proof of the first part of Theorem 4.4
needs the theory of lifting property in [27, 28], and hence also the
completeness of $(\Omega,{\cal F},P)$. The second part of it is most
important and need not assume that $(\Omega,{\cal F},P)$ is
complete, since $L^{0}({\cal F},B)$ and $L^{0}(\hat{{\cal F}},B)$ as
well as $L^{0}({\cal F},B')$ and $L^{0}(\hat{{\cal F}},B')$ can be
identified, so that we can first prove that $L^{0}(\hat{{\cal
F}},B')\cong(L^{0}(\hat{{\cal F}},B))^{\ast}$
iff $B^{\prime}$ has the Radon-Nikod\'{y}m property with respect to
$(\Omega, {\cal F}, P)$, and then return to our desired result.

To give Theorem 4.8 below, we first give Theorem 4.6. Let $1\leq
p\leq \infty$ and $(S, \|\cdot\|)$ an $RN$ module over $K$ with base
$(\Omega, {\cal F}, P)$. Define $\|\cdot\|_{p}:S\rightarrow
[0,+\infty]$ by $\|x\|_{p}=(\int_{\Omega}\|x\|^{p}dP)^{\frac{1}{p}}$
for $p\in [1,+\infty)$ and $\|\cdot\|_{\infty}=$ the essential
supremum of $\|x\|$, $\forall x\in S$, and denote by
$L^{p}(S)=\{x\in S\,|\, \|x\|_{p}<+\infty \}$, then $(L^{p}(S),
\|\cdot\|_{p})$ is a normed space, for all $q,1\leq q\leq+\infty$ we
can also have  $(L^{q}(S^{\ast}), \|\cdot\|_{q})$ in a similar way.
The following Theorem 4.6 is essentially independent of a special
choice of ${\cal T}_{c}$ and ${\cal T}_{\varepsilon,\lambda}$, which
was proved in [10] under ${\cal T}_{\varepsilon,\lambda}$, whose
proof in English was given in [14].

{\noindent\bf Theorem 4.6 [10, 14].} Let $1\leq p<+\infty$ and $1<
q\leq+\infty$ be a pair of H\"{o}lder conjugate numbers. Then
$(L^{q}(S^{\ast}), \|\cdot\|_{q})$ is isometrically isomorphic
with the classical
conjugate space of $(L^{p}(S), \|\cdot\|_{p})$, denoted by
$(L^{p}(S))^{\prime}$, under the canonical mapping $T:
L^{q}(S^{\ast})\rightarrow (L^{p}(S))^{\prime}$ defined as follows.
For each $f\in L^{q}(S^{\ast})$, $T_{f}$, denoting $T(f)$, :
$L^{p}(S)\rightarrow K$ is defined by $T_{f}(g)=\int_{\Omega}f(g)dP$
for all $g\in L^{p}(S)$.

Theorem 4.6 gives all representation theorems of the dual of
Lebesgue-Bochner function spaces by taking $S=L^{0}({\cal F}, B)$
(at which time $L^{p}(S)$ is just $L^{p}({\cal F}, B)$, the
classical Lebesgue-Bochner function spaces, see [12] for more
details). On the other hand, it provides the connection between the
random conjugate space $S^{\ast}$ and the classical conjugate space
$(L^{p}(S))^{\prime}$, and thus a powerful tool for the theory of
random conjugate spaces, cf.[10, 19, 22]. In this paper, we will
still employ it in the proof of Theorem 4.8. In particular,
combining the ideas of constructing
$L^{p}_{{\cal F}}(\mathcal {E})$ in [4,25] and $L^{p}(S)$ as above
at once leads us to the following:

{\noindent\bf Example 4.7.} Let $(S, \|\cdot\|)$ be an $RN$ module
over $K$ with base $(\Omega, \mathcal{E}, P)$ and ${\cal F}$ a
sub-$\sigma-$algebra of $\mathcal{E}$. For each $1\leq
p\leq+\infty$, define $|||\cdot|||_{p}: S\rightarrow L_{+}^{0}({\cal
F}, \bar{R})=\{\xi\in L^{0}({\cal F}, \bar{R}) ~|~ \xi\geq
0\}$ as follows: for all $x\in S$,

\[ |||x|||_{p}=\left\{
   \begin{array}{ll}
   (E[\|x\|^{p}\,|\,{\cal F}])^{\frac{1}{p}},    &\mbox{$if$ $p\in [1, +\infty)$}\\
    \wedge\{\xi\in L_{+}^{0}({\cal F}, \bar{R})\,|\, \xi\geq
    \|x\|\},    &\mbox{$if$ $p=+\infty$}
   \end{array}
   \right.,
\]\\
where $E[\cdot\,|\,{\cal F}]$ denotes the conditional expectation,
cf. [4,25].

Denote $L^{p}_{{\cal F}}(S)=\{x\in E\,|\, |||x|||_{p}\in
L_{+}^{0}({\cal F}, R)\}$, then $(L^{p}_{{\cal F}}(S),
|||\cdot|||_{p})$ is an $RN$ module over $K$ with base
$(\Omega,{\cal F},P)$. Similarly, we can also have $L^{q}_{{\cal
F}}(S^{\ast})$ for all $q\in [1, +\infty]$. When
$S=L^{0}(\mathscr{E}, R)$, $L^{p}_{{\cal F}}(S)$ is exactly
$L^{p}_{{\cal F}}(\mathscr{E})$ of [4, 25].

{\noindent\bf Theorem 4.8.} Let $1\leq p<+\infty$ and $1<
q\leq+\infty$ be a pair of H\"{o}lder conjugate numbers. The
canonical mapping $T:L^{q}_{{\cal F}}(S^{\ast})\rightarrow
(L^{p}_{{\cal F}}(S))^{\prime}$ is surjective and random-norm
preserving
, where for each $f\in L^{q}_{{\cal F}}(S^{\ast})$,
$T_{f}$(denoting $T(f)$): $L^{p}_{{\cal F}}(S)\rightarrow
L^{0}({\cal F}, K)$ is defined by $T_{f}(g)=E[f(g)\,|\,{\cal F}]$
for all $g\in L^{p}_{{\cal F}}(S)$, and $L^{q}_{{\cal
F}}(S^{\ast})$ and $L^{p}_{{\cal F}}(S)$ are the same as in
Example 4.7.

For the sake of clearness, the proof of Theorem 4.8 is divided into
the following two Lemmas 4.9 and 4.10, Lemma 4.9 shows that $T$ is
well defined and isometric (namely random-norm preserving) and Lemma
4.10 shows that $T$ is surjective. Specially, we need to remind the
readers of noticing that $|||\cdot|||_{p}$ and $|||\cdot|||_{q}$ are
the $L^{0}-$norms on $L^{p}_{{\cal F}}(S)$ and on $L^{q}_{{\cal
F}}(S^{\ast})$, respectively, whereas $\|\cdot\|_{p}$ and
$\|\cdot\|_{q}$ are norms.

{\noindent\bf Lemma 4.9.} T is well defined and isometric.

{\noindent\bf Proof.} For any fixed $f\in L^{q}_{{\cal
F}}(S^{\ast})$, we will
first prove that $T_f\in(L^{p}_{{\cal F}}(S))^{\ast}$ and
$\|T_f\|=||| f|||_q$ when $p>1$ as follows:

For any $g\in L^{p}_{{\cal F}}(S), T_f(g)=E[f(g)\,|\,{\cal F}]$,
we have the following:
\begin{eqnarray*}
|T_f(g)|&\leq& E[|f(g)|\,|\,{\cal F}]\\
&\leq& E[\|f\|\cdot\|g\|\,|\,{\cal F}]\\
&\leq& |||
f|||_q\cdot(E[\|g\|^{p}\,|\,{\cal F}])^{\frac{1}{p}} \hspace{8 cm}(4.26)\\
&=& ||| f|||_q\cdot||| g|||_p
\end{eqnarray*}

This shows that $T_f\in(L^{p}_{{\cal F}}(S))^{\ast}$ and
$\|T_f\|\leq||| f|||_q$, namely $T$ is well defined.
We remain to prove $\|T_f\|=||| f|||_q$ when $p>1$.

Let $\xi$ be an arbitrary representative of $|||
f|||_q$ and
$A_n=\{\omega\in\Omega~|~n-1\leq\xi(\omega)<n\}$ for each $n\in
N$. Then $\{A_n~|~n\in N\}$ forms a countable partition of
$\Omega$ to ${\cal F}$. Observing $\int_{\Omega}|||
g|||_p^{p}dP=\int_{\Omega}E[\|g\|^{p}\,|\,{\cal
F}]dP=\int_{\Omega}\|g\|^{p}dP,\forall g\in L^{p}_{{\cal F}}(S)$,
thus we have the following relation:

$L^{p}(L^{p}_{{\cal F}}(S))=L^{p}(S)\hfill{}(4.27)$

Since $p>1, 1<q<+\infty,$ we also have the relation:

$L^{q}(L^{q}_{{\cal F}}(S^{\ast}))=L^{q}(S^{\ast})\hfill{}(4.28)$

Now, we fix $n$ and prove

$\tilde{I}_{A_n}\|T_f\|=\tilde{I}_{A_n}|||
f|||_q\hfill{}(4.29)$

Since $\tilde{I}_{A_n}T_f(g)=E[\tilde{I}_{A_n}f(g)\,|\,{\cal F}]$
for all $g\in L^{p}_{{\cal F}}(S)$, we have, of course, that
$\tilde{I}_{A_n}T_f(g)=E[\tilde{I}_{A_n}f(g)\,|\,{\cal F}]$ for
all $g\in L^{p}(L^{p}_{{\cal F}}(S))=L^{p}(S)$. Further, since
$\tilde{I}_A L^{p}(L^{p}_{{\cal F}}(S))=\tilde{I}_A
L^{p}(S)\subset L^{p}(L^{p}_{{\cal F}}(S))=L^{p}(S)$ for all $A\in
{\cal F}$, we can have the following important relation:

$\tilde{I}_{A_n}\tilde{I}_{A}T_f(g)=E[\tilde{I}_{A_n}\tilde{I}_{A}f(g)\,|\,{\cal F}],\hfill{}(4.30)$\\
for all $A\in{\cal F}$ and all $g\in L^{p}(L^{p}_{{\cal
F}}(S))=L^{p}(S)$

Obviously, from (4.30), for all $A\in{\cal F}$ we can have the
following relation:

$\int_{\Omega}(\tilde{I}_{A_n}\tilde{I}_{A}T_f)(g)dP=\int_{\Omega}(\tilde{I}_{A_n}\tilde{I}_{A}f)(g)dP,\hfill{}(4.31)$\\
for all $g\in L^{p}(L^{p}_{{\cal F}}(S))=L^{p}(S).$

For each fixed $A\in{\cal F}$, the left side of (4.31) defines a
bounded linear functional on $L^{p}(L^{p}_{{\cal F}}(S))$, whose
norm is equal to
$(\int_{\Omega}\|\tilde{I}_{A_n}\tilde{I}_{A}T_f\|^{q}dP)^{\frac{1}{q}}=(\int_{A}(\tilde{I}_{A_n}\|T_f\|)^{q}dP)^{\frac{1}{q}}$
by applying Theorem 4.6 to $L^{p}(L^{p}_{{\cal F}}(S))$. The same
bounded linear functional is also a bounded linear functional on
$L^{p}(S)$ defined by the right side of (4.31), then whose norm is
also equal to
$(\int_{A}(\tilde{I}_{A_n}\|f\|)^{q}dP)^{\frac{1}{q}}$.

Consequently,
$\int_{A}\|\tilde{I}_{A_n}T_f\|^{q}dP=\int_{A}(\tilde{I}_{A_n}\|f\|)^{q}dP$
for all $A\in{\cal F}$. Since $\|\tilde{I}_{A_n}T_f\|^{q}\in
L^{0}_{+}({\cal F},R)$, we have
$\|\tilde{I}_{A_n}T_f\|^{q}=E[(\tilde{I}_{A_n}\|f\|)^{q}\,|\,{\cal
F}]$. Again noticing $A_n\in{\cal F}$, we can have
$\tilde{I}_{A_n}\|T_f\|=\tilde{I}_{A_n}(E[\|f\|^{q}\,|\,{\cal
F}])^{\frac{1}{q}}$, which is just (4.29).

Since $\sum_{n\geq 1}A_n=\Omega, \|T_f\|=|||
f|||_q$.

Finally, we consider the case of $p=1$, for the sake of clearness,
we use $|\cdot|_{\infty}$ for the usual $L^{\infty}-$norm on the
Banach space $L^{\infty}({\cal E},K)$ of equivalence classes of
essentially bounded ${\cal E}-$measurable $K-$valued functions on
$(\Omega,{\cal E},P)$. Then it is easy to see that
$L^{\infty}(L^{\infty}_{{\cal F}}(S^{\ast}))=L^{\infty}(S^{\ast})$,
in particular that $\|f\|_{\infty}=||||
f|||_{\infty}|_{\infty}$ for all $f\in
L^{\infty}(S^{\ast})$.

Since, we can, similarly to the case $p>1$, have , by noticing
$\tilde{I}_{A_n}f\in L^{\infty}(S^{\ast}),$ the following relation:
$|\tilde{I}_A(\tilde{I}_{A_n}\|T_f\|)|_{\infty}=|(\tilde{I}_{A_n}\tilde{I}_A\|f\|)|_{\infty}$
for all $A\in {\cal F}$, but as stated above, the latter is just
equal to
$|\tilde{I}_{A}|||\tilde{I}_{A_n}f|||_{\infty}|_{\infty}$.
Since $\tilde{I}_{A_n}\|T_f\|$ and
$|||\tilde{I}_{A_n}f|||_{\infty}$(namely
$\tilde{I}_{A_n}||| f|||_{\infty}$) are both in
$L^{0}_{+}({\cal F},R)$, then we must have that
$\tilde{I}_{A_n}\|T_f\|=\tilde{I}_{A_n}|||
f|||_{\infty}$, again since $\sum_{n\geq 1}A_n=\Omega,$ we
can have $\|T_f\|=||| f|||_{\infty}.$ \quad $\Box$

{\noindent\bf Lemma 4.10.} $T$ is surjective.

{\noindent\bf Proof.} Let $F$ be an arbitrary element of
$(L^{p}_{{\cal F}}(S))^{\ast}$. We want to prove that there exists
an $f\in L^{q}_{{\cal F}}(S^{\ast})$ such that $F=T_f$.

Since $\|F\|\in L^{0}_{+}({\cal F},R)$, letting $\xi$ be a chosen
representative of $\|F\|$ and
$A_n=\{\omega\in\Omega\,|\,n-1\leq\xi(\omega)<n\}$ for each $n\in
N$, then $\{A_n\,|\,n\in N\}$ forms a countable partition of
$\Omega$ to ${\cal F}$. Since $|F(g)|\leq\|F\||||
g|||_p=\|F\|(E[\|g\|^{p}\,|\,{\cal F}])^{\frac{1}{p}},
|\tilde{I}_{A_n}F(g)|\leq n(E[\|g\|^{p}\,|\,{\cal
F}])^{\frac{1}{p}},\forall g\in L^{p}_{{\cal F}}(S)$ and $n\in N$.

Now, we fix $n$ and notice
$|\int_{\Omega}(\tilde{I}_{A_n}F)(g)dP|\leq
n\int_{\Omega}(E[\|g\|^{p}\,|\,{\cal F}])^{\frac{1}{p}}dP\leq
n(\int_{\Omega}\|g\|^{p}dP)^{\frac{1}{p}}$ for all $g\in L^p(S)$
by the  H\"{o}lder inequality, then by Theorem 4.6 there exists
$f_{n}\in L^{q}(S^{\ast})$  such that$
\int_{\Omega}(\tilde{I}_{A_n}F)(g)dP=\int_{\Omega}f_{n}(g)dP$ for
all $g\in L^{p}(S)$.

Since $\tilde{I}_{A}L^{p}(S)\subset L^{p}(S)$ for all $A\in {\cal
F}$, we have that
$\int_{A}(\tilde{I}_{A_n}F)(g)dP=\int_{\Omega}(\tilde{I}_{A_n}F)(\tilde{I}_{A}g)dP=\int_{\Omega}f_{n}(\tilde{I}_{A}g)dP=\int_{A}f_{n}(g)dP$
for all $g\in L^p(S)$ and all $ A\in {\cal F}$, which yields the
following important relation (by noticing $\tilde{I}_{A_n}F(g)\in
L^{0}({\cal F},K)$):
$$\tilde{I}_{A_n}F(g)=E[f_n(g)\,|\,{\cal F}],\eqno(4.32)$$
for all $g\in L^{p}(S)\equiv L^{p}(L^{p}_{\cal F}(S))$.

Since $f_n\in L^{q}(S^\ast)\subset L^{q}_{\cal F}(S^\ast)$, then
Lemma 4.9 shows that $T_{f_{n}}\in(L^{q}_{\cal F}(S))^{\ast}$, and
(4.32) just shows that $\tilde{I}_{A_n}F$ and $T_{f_{n}}$ are equal
on $L^{p}(L^{p}_{\cal F}(S))$. Since $\tilde{I}_{A_n}F$ and
$T_{f_{n}}$ are both in $(L^{p}_{\cal F}(S))^{\ast}=(L^{p}_{\cal
F}(S))^{\ast}_{\varepsilon,\lambda}$, namely they are both
continuous module homomorphism from $(L^{p}_{{\cal F}}(S),{\cal
T}_{\varepsilon,\lambda})$ to $(L^{0}({\cal F},K),{\cal
T}_{\varepsilon,\lambda})$, and
$L^{p}(L^{p}_{{\cal F}}(S))$ is ${\cal
T}_{\varepsilon,\lambda}-$dense in $L^{p}_{{\cal F}}(S)$
(cf.[19,22]), $\tilde{I}_{A_n}F=T_{f_{n}}$, namely we can have the
following relation:
                       $$\tilde{I}_{A_n}F(g)=E[f_{n}(g)\mid{\cal F}],\eqno(4.33)$$
                       for all $ g\in L^{p}_{{\cal F}}(S)$.

We have from (4.33) the following relation:
$$\tilde{I}_{A_n}F(g)=E[\tilde{I}_{A_n}f_{n}(g)\mid{\cal F}],\eqno(4.34)$$ for all $ g\in L^{p}_{{\cal F}}(S)$.

Let $ f=\sum_{n\geq1}\tilde{I}_{A_n}f_{n}$. Since $L^{q}_{{\cal
F}}(S^{\ast})$ has the countable concatenation property ( or we can
directly define $ f(g)=\sum_{n\geq1}\tilde{I}_{A_n}f_{n}(g),\forall
g\in L^{p}_{{\cal F}}(S)$ and verify that $f$ is first in $S^{\ast}$
and then $f\in L^{q}_{{\cal F}}(S^{\ast})$), we have that $f\in
L^{q}_{{\cal F}}(S^{\ast})$.

Finally, (4.34) shows that $F(g)=E[f(g)\mid{\cal
F}]=T_{f}(g),\forall g\in L^{p}_{{\cal F}}(S)$. $\Box$

{\noindent \bf{Remark 4.11.}} Let $K^{d}$ be the $d-$dimensional
Euclidean space over $K$. Then Theorem 4.3 and the second part of
Theorem 4.4 both imply that $(L^{0}({\cal F},K^{d}))^{\ast}$
=$L^{0}({\cal F},K^{d})$, which is just Proposition 4.2 of [25].
Since $(L^{0}({\cal E},K))^{\ast}=L^{0}({\cal E},K)$, if we take
$S=L^{0}({\cal E},K)$ in Theorem 4.8, then we have $(L^{p}_{{\cal
F}}({\cal E}))^{\ast}=L^{q}_{{\cal F}}({\cal E})$, which is just
Theorem 4.5 of [25], so our Theorem 4.8 is surprisingly general.
Besides, the proof of the isometric property of $T$ of Theorem 4.8
is completely new since Theorem 4.5 of [25] did not involve any
isometric arguments. Further, we can also have $L^{p}_{{\cal
F}}(S)=L^{0}({\cal F},K)\cdot L^{p}(S)$.

\subsection*{4.2. Random reflexivity and the James theorem}

\quad\, Let $(E,\|\cdot\|)$ be an $RN$ module, $E^{\ast\ast}$
denotes $(E^{\ast})^{\ast}$, the canonical embedding mapping
$J:E\rightarrow E^{\ast\ast}$ defined by $(Jx)(f)=f(x),\forall x\in
E$ and $f\in E^{\ast}$, is random$-$norm preserving. If $J$ is
surjective, then $E$ is called {\it{random reflexive}}.

Clearly, if $E$ is random reflexive, then $E$ is complete under both
${\cal T}_c$ and ${\cal T}_{\varepsilon,\lambda}$, and has the
countable concatenation property since $E^{\ast\ast}$ has these
properties, and the main results of this subsection are essentially
independent of a special choice of ${\cal T}_c$ and ${\cal
T}_{\varepsilon,\lambda}$, which were established under ${\cal
T}_{\varepsilon,\lambda}$. Since they are still valid under ${\cal
T}_c$, we only state them without proofs except Theorem 4.13 and
without mention of topologies.

{\noindent\bf Theorem 4.12 ([8]).} $L^{0}({\cal F},B)$ (see Theorem
4.4) is random reflexive iff $B$ is a reflexive Banach space.

{\noindent\bf Theorem 4.13.} Let $(S,\|\cdot\|)$ be an $RN$ module
over $K$ with base $(\Omega,{\cal E},P),1<p<+\infty$ and ${\cal F}$
a sub-$\sigma-$algebra of ${\cal E}$. Then the following statements
are equivalent:

(1) $(S,\|\cdot\|)$ is random reflexive;

(2) $L^{p}(S)$ is a reflexive Banach space;

(3) $L^{p}_{{\cal F}}(S)$ is random reflexive.

{\noindent\bf Proof.} $(1)\Leftrightarrow (2)$ was proved in [10],
see also [14] for its proof in English. $(2) \Leftrightarrow (3)$
has been implied by $(1)\Leftrightarrow (2)$ by noticing
$L^{p}(L^{p}_{{\cal F}}(S))=L^{p}(S)!$ \quad$\Box$

{\noindent\bf Theorem 4.14 [19]}. A complete $RN$ module
$(S,\|\cdot\|)$ is random reflexive iff their exists $p\in S(1)$ for
each $f\in E^{\ast}$ such that $f(p)=\|f\|$, where $S(1)=\{p\in
S\,|\,\|p\|\leq 1\}$.

\subsection*{4.3. Banach-Alaoglu Theorem and
 Banach-Bourbaki-Kakutani-$\check{\textmd{S}}$mulian Theorem}
\quad\, Let $(E,\|\cdot\|)$ be an $RN$ module over $K$ with base
$(\Omega,{\cal F},P)$. Let $\xi=\vee\{\|x\|\,|\,x\in E\}$ and let
$\xi^{0}$ be a representative of $\xi$, the set
$\{\omega\in\Omega\,|\,\xi^{0}(\omega)>0\}$ is called the support of
$E$ (unique a.s.), if $\Omega$ is the support, then $E$ is called
{\it{having full support}}, in this subsection $(E,\|\cdot\|)$ is
always assumed to have full support. $A\in{\cal F}$ is called a
$P-$atom if $P(A)>0$, and if $B\in{\cal F}$ and $B\subset A$ must
imply either $P(B)=0$ or $P(A\backslash B)=0$. $(\Omega,{\cal F},P)$
is said to be essentially purely $P-$atomic if there is a sequence
$\{A_n\,|\,n\in N\}$ of disjoint $P-$atoms such that $\sum_{n\geq
1}A_n=\Omega$ and ${\cal F}\subset\overline{\sigma\{A_n\,|\,n\in
N\}}^{P}$, where $\overline{\sigma\{A_n\,|\,n\in N\}}^{P}$ denotes
the $P-$completion of the $\sigma-$algebra generated by
$\{A_n\,|\,n\in N\}$.

The two new results of the subsection are Theorems 4.16 and 4.18,
and Theorems 4.15 and 4.17 are stated in order to contrast with the
former two.

For the four topologies
$\sigma_{\varepsilon,\lambda}(E,E^{\ast}),\sigma_c(E,E^{\ast}),\sigma_{\varepsilon,\lambda}(E^{\ast},E)$
and $\sigma_c(E^{\ast},E)$, we refer to Definition 3.25.

Classical Banach-Alaoglu theorem says the closed unit ball
$B^{\prime}(1)$ of the conjugate space $B^{\prime}$ of a normed
space $B$ is always $\sigma(B^{\prime},B)-$compact, namely
$\mbox{w}^{\ast}-$compact, but for an $RN$ module $(E,\|\cdot\|)$
with base $(\Omega,{\cal F},P)$ we have the following:

{\noindent\bf Theorem 4.15 ([16]).} $E^{\ast}(1)=\{f\in
E^{\ast}\,|\,\|f\|\leq 1\}$ is
$\sigma_{\varepsilon,\lambda}(E^{\ast},E)-$compact iff ${\cal F}$ is
essentially
purely $P-$atomic.

{\noindent\bf Theorem 4.16}. If $E^{\ast}(1)$ is
$\sigma_{c}(E^{\ast},E)-$compact then ${\cal F}$ is essentially
purely $P-$atomic. But the converse is not true.

{\noindent\bf Proof.} Since $\sigma_{c}(E^{\ast},E)$ is stronger
than $\sigma_{\varepsilon,\lambda}(E^{\ast}, E), E^{\ast}(1)$ is
$\sigma_{\varepsilon,\lambda}(E^{\ast},E)-$compact, namely ${\cal
F}$ is essentially purely $P-$atomic if $E^{\ast}(1)$ is
$\sigma_{c}(E^{\ast},E)-$compact.

From the process of proof of Theorem 4.15 given in [16] we can
similarly prove that $E^{\ast}(1)$ is
$\sigma_{c}(E^{\ast},E)-$compact iff $\{x\in L^{0}({\cal
F},K)\,|\,|x|\leq1\}$ is ${\cal T}_{c}-$compact. We can construct
the following example to show that even if $(\Omega,{\cal F},P)$ is
essentially purely $P-$atomic $\{x\in L^{0}({\cal
F},K)\,|\,|x|\leq1\}$ is not ${\cal T}_{c}-$compact, either. Take
$\Omega=N,{\cal F}=2^{N}$ and $P(A)=\sum_{i\in A}\frac{1}{2^{i}}$
for all $A\in {\cal F}$, then $(\Omega,{\cal F},P)$ is purely
atomic, but at this time $\{x\in L^{0}({\cal F},K)\,|\,|x|\leq1\}$
is exactly the closed unit ball of the Banach space $l^{\infty}$ of
bounded sequences in $K$, it is well known that it is not
norm$-$compact in the Banach space, and hence it is not ${\cal
T}_{c}-$compact, either, since ${\cal T}_{c}-$topology is stronger than the
norm$-$topology on the closed unit ball, so that $E^{\ast}(1)$ is
not $\sigma_{c}(E^{\ast},E)-$compact. $\Box$

Classical Banach-Bourbaki-Kakutani-$\check{\textmd{S}}$mulian
theorem says that a Banach space is reflexive iff its closed unit
ball is weakly compact, in [16] we prove the following:

{\noindent\bf Theorem 4.17([16])}. $E(1)=\{x\in E\,|\,\|x\|\leq1\}$
is $\sigma_{\varepsilon,\lambda}(E,E^{\ast})-$compact iff both
$(E,\|\cdot\|)$ is random reflexive and ${\cal F}$ is essentially
purely $P-$atomic.

Just as the norm$-$topology and weak one on $K^{d}$ are the same, it
is easy to check that ${\cal
T}_{\varepsilon,\lambda}=\sigma_{\varepsilon,\lambda}(L^{0}({\cal
F},K^{d}),L^{0}({\cal F},K^{d}))$ and ${\cal
T}_c=\sigma_c(L^{0}({\cal F},K^{d}),L^{0}({\cal F},K^{d}))$ on
$L^{0}({\cal F},K^{d})$. From this it is easy to see the proof of
the second part of Theorem 4.18 below.

{\noindent\bf Theorem 4.18.} If $E(1)$ is
$\sigma_{c}(E,E^{\ast})-$compact then $(E,\|\cdot\|)$ is random
reflexive and ${\cal F}$ is essentially purely $P-$atomic. But the
converse is not true.

{\noindent\bf Proof.} The proof is similar to the one of Theorem
4.16., in particular, if we take $E=L^{0}({\cal F},K)$, then that
counterexample in the proof of Theorem 4.6 can also serve for the
proof of the converse of the theorem. \quad$\Box$

\section*{5. Some further remarks on the $(\varepsilon,\lambda)-$topology
and the locally $L^{0}-$convex topology}

\quad\, The topology of convergence in probability is one of the
most useful topologies on the space of random variables. The
$(\varepsilon,\lambda)-$topology, as a natural generalization of
the former, makes the theory of $RN$ and $RIP$ modules naturally
applicable to many topics in probability theory, for example, our
recent work [18] provides some interesting applications of $RIP$
modules to random linear operators on Hilbert spaces (cf.[27]).
Further it has many advantages itself. For example, it admits a
countable concatenation skill in ${\cal
T}_{\varepsilon,\lambda}-$complete random locally convex modules
(see Example 3.7), and has many nice properties, for example,
$L^{p}(S)$ is ${\cal T}_{\varepsilon,\lambda}-$dense in $S$ for an
$RN$ module $S$, which produces the useful connection of random
conjugate spaces with classical conjugate spaces [see Theorem 4.6
and the process of proof of Theorem 4.8]. The
$(\varepsilon,\lambda)-$topology is in harmony with the module
structure, the family of $L^{0}-$seminorms and the order structure
on $L^{0}({\cal F},R)$ of a random locally convex module so that a
random locally convex module and its random conjugate space can be
deeply developed under the framework of topological modules over
the topological algebra $(L^{0}({\cal F},K),{\cal
T}_{\varepsilon,\lambda})$. However, the
$(\varepsilon,\lambda)-$topology on the linear spaces is rarely a
locally convex topology, which makes us not establish such results
as Theorem 2.6 and Lemma 3.10 of [4], in particular the
$(\varepsilon,\lambda)-$topology can neither perfectly match the
notion of a gauge function introduced in [4].

The locally $L^{0}-$convex topology has the nice $L^{0}-$convexity
and perfectly matches the gauge functions, which admits Theorem 2.6
and Lemma 3.10 of [4] and thus has played a crucial role in convex analysis
on $L^{0}-$modules. We can predict that the locally $L^{0}-$convex
topology will also develop its power in non-$L^{0}-$linear analysis.
However, it is too
strong to make the previous frequently used skills reserved, for
example, the countable concatenation skill often fails, it is
impossible that $L^{p}(S)$ is ${\cal T}_c-$dense in $S$, and in
particular it is also impossible to establish Theorem 3.12 under
${\cal T}_c$.

Comprehensively speaking, the two kinds of topologies can be both
applied to mathematical finance, cf.[1,4] and the references
therein, and the principal results of this paper enough convince
people that the two kinds of topologies should be, simultaneously
rather than in a single way, considered in the future study of
random locally convex modules together with their financial
applications.
\section*{Acknowledgements}
\quad\, The author is greatly indebted to Professor Berthold
Schweizer for ten years of invaluable suggestions and support to
the author's work. He is also greatly indebted to Professor Daniel
W. Stroock for giving the author some invaluable directions to
mathematics and generous financial support during a short stay at
MIT from November of 2001 to March of 2002. He would like to
sincerely thank Professors Shijian Yan, Yingming Liu, Kung-Ching
Chang, Shunhua Sun, Mingzhu Yang, Peide Liu, Fuzhou Gong and Kunyu
Guo for their advices and kindness given on some occasions.
Finally, the author would like to thank his students Xiaolin Zeng,
Xia Zhang, Shien Zhao and Yujie Yang for some helps in the process
of preparing this paper.

\vskip0.2in
\no {\bf References}
\vskip0.1in

\footnotesize
\REF{[1] } W.Brannath, W.Schachermayer, A bipolar theorem for
subsets of $L^{0}_{+}(\Omega,{\cal F},P)$, in: S\'{e}minaire de
Probabilit\'{e}s \uppercase\expandafter{\romannumeral 33}, in:
Lecture Notes in Math., vol. 1709, Springer, 1999, pp.349--354.

\REF{[2] } W.W.Breckner, E.Scheiber, A Hahn-Banach extension theorem
for linear mappings into ordered modules, Mathematica 19(42)(1977)
13--27.

\REF{[3] } N. Dunford, J. T. Schwartz, Linear Operators(\uppercase
\expandafter {\romannumeral 1}), Interscience, New York, 1957.

\REF{[4] } D.Filipovi\'{c}, M.Kupper, N.Vogelpoth, Separation and
duality in locally $L^{0}-$convex modules, J. Funct. Anal. 256(2009)
3996--4029.

\REF{[5] } T.X. Guo, The theory of probabilistic metric spaces with
applications to random functional analysis, Master's thesis, Xi'an
Jiaotong University (China), 1989.

\REF{[6] } T.X. Guo, Random metric theory and its applications, Ph.D
\ Thesis, Xi'an Jiaotong University (China), 1992.

\REF{[7] } T.X. Guo, Extension theorems of continuous random linear
operators on random domains, J. Math. Anal. Appl. 193(1)(1995)
15--27.

\REF{[8] } T.X. Guo, The Radon-Nikod$\acute{y}$m property of
conjugate spaces and the $\mbox{w}^{\ast}-$equivalence theorem for
$\mbox{w}^{\ast}-$ measurable functions, Sci. China Ser. A-Math.
39(1996) 1034--1041.

\REF{[9] } T.X. Guo, Module homomorphisms on random normed modules,
Chinese Northeastern Math. J. 12(1996) 102--114.

\REF{[10] } T.X. Guo, A characterization for a complete random
normed module to be random reflexive, J. Xiamen Univ. Natur. Sci.
36(1997) 499--502 (In Chinese).

\REF{[11] } T.X. Guo, Some basic theories of random normed linear
spaces and random inner product spaces, Acta Anal. Funct.
 Appl. 1(2)(1999) 160--184.

\REF{[12] } T.X. Guo, Representation theorems of the dual of
Lebesgue-Bochner function spaces, Sci. China Ser. A-Math. 43(2000)
234--243.

\REF{[13] } T.X. Guo, Survey of recent developments of random metric
theory and its applications in China(I), Acta Anal. Funct. Appl.
3(2001) 129--158.

\REF{[14] } T.X. Guo, Survey of recent developments of random metric
theory and its applications in China(II), Acta Anal. Funct. Appl.
3(2001) 208--230.

\REF{[15] } T.X. Guo, Several applications of the theory of random
conjugate spaces to measurability problems, Sci. China Ser. A-Math.
50(2007) 737--747.

\REF{[16] } T.X. Guo, The relation of Banach-Alaoglu theorem and
Banach-Bourbaki-Kakutani-$\check{\textmd{S}}$mulian theorem in
complete random normed modules to stratification structure, Sci
China(Ser A)51(2008) 1651--1663.

\REF{[17] } T.X. Guo, X.X. Chen, Random duality, Sci. China Ser.
A(2009), in press.

\REF{[18] } T.X. Guo, F.Z. Gong, the random spectral representation
theorem of random self-adjoint operators on complete random inner
product modules, to appear.

\REF{[19] } T.X. Guo, and S.B. Li, The James theorem in complete
random normed modules, J. Math. Anal. Appl. 308(2005) 257--265.

\REF{[20] } T.X. Guo, S.L. Peng, A characterization for an
$L(\mu,K)-$topological module to admit enough canonical module
homomorphisms, J. Math. Anal. Appl. 263(2001) 580--599.

\REF{[21] } T.X. Guo, H.X. Xiao, A separation theorem in random
normed modules, J. Xiamen Univ. Natur. Sci. 42(2003) 270--274.

\REF{[22] } T.X. Guo, H.X. Xiao, X.X. Chen, A basic strict
separation theorem in random locally convex modules, Nonlinear Anal.
71(2009) 3794--3804.

\REF{[23] } T.X. Guo, Z.Y. You, The Riesz's representation theorem
in complete random inner product modules and its applications,
Chinese. Ann. of Math. Ser. A. 17(1996) 361--364.

\REF{[24] } T.X. Guo, L.H. Zhu, A characterization of continuous
module homomorphisms on random seminormed modules and its
applications, Acta Math. Sinica English Ser. 19(1)(2003) 201--208.

\REF{[25] } M.Kupper, N.Vogelpoth, Complete $L^{0}-$normed modules
and automatic continuity of monotone convex functions, Working paper
Series No.10, Vienna Institute of Finance, 2008.

\REF{[26] } B.Schweizer, A.Sklar, Probabilistic Metric Spaces,
Elsevier/North$-$Holland, New York, 1983; reissued by Dover Publications, Mineola, New York, 2005.

\REF{[27] } A.V.Skorohod, Random Linear Operators, D.Reidel
Publishing Company, Holland, 1984.

\REF{[28] } A.Ionescu Tulcea, C.Ionescu Tulcea, On the lifting
property(\uppercase\expandafter{\romannumeral 1}). J. Math. Anal.
Appl. 3(1961) 537-546.

\REF{[29] } A.Ionescu Tulcea, C.Ionescu Tulcea, On the lifting
property(\uppercase\expandafter{\romannumeral 2}). J. Math. Mech.
11(5)(1962) 773-795.

\REF{[30] } D.Vuza, The Hahn-Banach theorem for modules over ordered
rings, Rev. Roumaine Math. Pures Appl. 9(27)(1982) 989--995.

\REF{[31] } D.H.Wagner, Survey of measurable selection theorems,
SIAM J. Control and Optimization, 15(1977) 859-903.

\end{document}